\documentclass[11pt]{article} 

\usepackage[utf8]{inputenc} 

\usepackage{geometry} 
\usepackage{graphicx} 

\usepackage{natbib}
\usepackage[most]{tcolorbox}
\usepackage{todonotes}
\usepackage{amssymb}
\usepackage{multirow}
\usepackage{algorithm}
\usepackage[noend]{algpseudocode}
\usepackage{amsmath}
\usepackage{array} 
\usepackage{bm}
\usepackage{caption}
\usepackage{subcaption}
\usepackage{hyperref}
\usepackage{glossaries}
\usepackage{booktabs} 
\usepackage{multicol}
\usepackage{paralist} 
\usepackage{setspace}
\usepackage{tikz}
\usetikzlibrary[shapes.geometric]
\usetikzlibrary{math}
\usepackage{optidef}
\usepackage{verbatim} 
\usepackage{ifthen}

\usepackage{fancyhdr} 
\usepackage{sectsty}
\allsectionsfont{\sffamily\mdseries\upshape} 

\usepackage[nottoc,notlof,notlot]{tocbibind} 
\usepackage[titles,subfigure]{tocloft} 

\newacronym{adg}{ADG}{action dependency graph}
\newacronym{coa}{CoA}{courses of action}
\newacronym{dag}{DAG}{directed acyclic graph}
\newacronym{dp}{DP}{dynamic programming}
\newacronym{lo}{LO}{linear optimization}
\newacronym{mdp}{MDP}{Markov decision process}
\newacronym{mio}{MIO}{mixed-integer optimization}
\newacronym{milo}{MILO}{mixed-integer linear optimization}
\newacronym{dt}{DT}{decision tree}
\newacronym{rl}{RL}{reinforcement learning}

\newacronym{fg}{FG}{full graph}
\newacronym{rg}{RG}{reduced graph}

\newcommand{\be}{\mathbf{e}}

\newcommand{\bs}{\mathbf{s}}
\newcommand{\bz}{\mathbf{z}}

\newcommand{\by}{\mathbf{y}}

\newcommand{\includeParagraph}{1} 

\usepackage{authblk}

\title{Efficient Tree Generation for Globally Optimal Decisions under Probabilistic Outcomes \thanks {\textbf{Approved for Public Release; Distribution Unlimited. Public Release Case Number 25-0045. \textcopyright 2025 The MITRE Corporation. ALL RIGHTS RESERVED. This technical data deliverable was developed using contract funds under Basic Contract No. W56KGU-18-D-0004.} The view, opinions, and/or findings contained in this report are those of The MITRE Corporation and should not be construed as an official Government position, policy, or decision, unless designated by other documentation."}}

\author{Berk \"Ozt\"urk\thanks{1ozturkbe@gmail.com}, She'ifa Punla-Green\thanks{spunla-green@mitre.org}, Les Servi\thanks{lservi@mitre.org} \\ The MITRE Corporation}

\date{April 10, 2025}

\begin{document}

\maketitle
\section*{Abstract}

Many real-world problems require making sequences of decisions where the outcomes of each decision are probabilistic and uncertain, and the availability of different actions is constrained by the outcomes of previous actions. There is a need to generate policies that are adaptive to uncertainty, globally optimal, and yet scalable as the state space grows. In this paper, we propose the generation of optimal decision trees, which dictate which actions should be implemented in different outcome scenarios, while maximizing the expected reward of the strategy. Using a combination of dynamic programming and mixed-integer linear optimization, the proposed methods scale to problems with large but finite state spaces, using problem-specific information to prune away large subsets of the state space that do not yield progress towards rewards. We demonstrate that the presented approach is able to find the globally optimal decision tree in linear time with respect to the number states explored.

\ifthenelse{\includeParagraph=0}{  
\section*{Executive Summary}
\label{sec:executivesummary}
This paper introduces a novel approach leveraging Operations Research optimization methods to identify the optimal decision tree of actions from a given catalog of alternatives actions, each characterized by probabilistic outcomes and potentially complex prerequisites dependent on the success or failure of prior actions. The methodology generates an actionable decision tree that dynamically guides users to the next optimal action based on the outcomes of previous actions. This framework incorporates budget constraints on actions and rewards tied to achieving specific end states, ensuring practical applicability in real-world scenarios.  If during the execution of the decision tree new information becomes available the approach permits the user to incorporate that new knowledge into the formulation and recompute the best future decision tree.

The proposed approach has broad applications across diverse domains, including healthcare, wargaming, and cybersecurity, where decision-making under uncertainty is found. Traditional tools such as Markov decision processes and reinforcement learning agents have shown limitations in handling intricate interdependencies between actions and distant reward horizons; this methodology overcomes these challenges. Specifically, they are designed to accommodate complex interactions and dependencies while maintaining computational efficiency.

A key advantage of this methodology is its scalability. The approach scales linearly with the number of states, enabling it to address very large systems that are often computationally prohibitive for other techniques. By combining rigorous optimization with practical scalability, this work offers a powerful and versatile solution for decision-making in complex, resource-constrained environments.
}

\ifthenelse{\includeParagraph=0}{  
\section*{Executive Summary}
\label{sec:executivesummary2}

Sequential decision making under uncertainty is a key characteristic of many decision problems in healthcare, wargaming and cyber. In these settings, an agent takes actions with discrete probabilistic outcomes that affect the state of its environment. Actions may or may not be available to the agent depending on the state of the environment, including the outcomes of previous actions. Given these complexities, the agent hopes to make progress towards rewards, attained via certain combinations of actions and associated outcomes. 

While tools like Markov decision processes and reinforcement learning agents have had some success in addressing certain classes of decision problems, they can have real drawbacks such as the inability to accommodate complex interdependencies and interactions between actions, their sensitivity to training parameters, and their frequent inability to learn good policies in face of distant reward horizons. In these settings, it is important to be able to have flexible tools that can handle a broad range of decision problems while guaranteeing globally optimal policies. 

The optimal strategy is conveyed to the user as a decision tree, where each node corresponds to a state and the optimal action to take next. Trees are a natural representation of possible alternative futures, and allow the user to easily 
}

\section{Introduction}
\label{sec:introduction}

Sequential decision making is a key aspect of many problems in domains such as wargaming, healthcare and cyber operations. In these settings, agents take actions to achieve goals, but the outcomes of each action are discrete, probabilistic and uncertain. This makes it challenging to know the best action in the face of alternative futures. Furthermore, available actions have complex conditional interdependencies that may constrain possible strategies, while leading to rewards in a distant time horizon. In this paper, we develop a method to define optimal decision sequences that adapt to outcomes of different actions, and represent them as decision trees. 

\begin{figure}[h!]
    \begin{center}
        \resizebox*{\textwidth}{!}{
        \begin{tikzpicture}
            \node[circle, line width=1mm, inner sep=0pt,
                minimum size = 6mm,
                draw=black!30] (root) at (-1.5,0) {$a_1$};
            \node[circle, line width=1mm, inner sep=0pt,
                minimum size = 6mm,
                draw=black!30] (l) at (-5,-2) {$a_3$};
            \node[circle, line width=1mm, inner sep=0pt,
                minimum size = 6mm,
                draw=black!30] (r) at (2.5,-2) {$a_4$};
            \node[circle, line width=1mm, inner sep=0pt,
                minimum size = 6mm,
                draw=black!30] (ll) at (-7,-4) {0};
            \node[circle, line width=1mm, inner sep=0pt,
                minimum size = 6mm,
                draw=black!30] (lr) at (-3,-4) {$a_7$};
            \node[circle, line width=1mm, inner sep=0pt,
                minimum size = 6mm,
                draw=black!30] (rl) at (0,-4) {$a_3$};
            \node[circle, line width=1mm, inner sep=0pt,
                minimum size = 6mm,
                draw=black!30] (rr) at (5,-4) {$a_5$};
            \node[circle, line width=1mm, inner sep=0pt,
                minimum size = 6mm,
                draw=black!30] (lrr) at (-2,-6) {100};
            \node[circle, line width=1mm, inner sep=0pt,
                minimum size = 6mm,
                draw=black!30] (rll) at (-1,-6) {0};
            \node[circle, line width=1mm, inner sep=0pt,
                minimum size = 6mm,
                draw=black!30] (lrl) at (-4,-6) {...};
            \node[circle, line width=1mm, inner sep=0pt,
                minimum size = 6mm,
                draw=black!30] (rlr) at (1,-6) {...};
            \node[circle, line width=1mm, inner sep=0pt,
                minimum size = 6mm,
                draw=black!30] (rrl) at (4,-6) {...};
            \node[circle, line width=1mm, inner sep=0pt,
                minimum size = 6mm,
                draw=black!30] (rrr) at (7,-6) {...};

            \draw[->, shorten >=1pt] (root) -- (l) node[draw=none,fill=none,font=\scriptsize, sloped,midway,below] {outcome 1};
            \draw[->, shorten >=1pt] (root) -- (r) node[draw=none,fill=none,font=\scriptsize, sloped,midway,below] {outcome 2};
            \draw[->, shorten >=1pt] (r) -- (rl) node[draw=none,fill=white,font=\scriptsize,midway] {1};
            \draw[->, shorten >=1pt] (r) -- (rr) node[draw=none,fill=white,font=\scriptsize,midway] {2};
            \draw[->, shorten >=1pt] (l) -- (lr) node[draw=none,fill=white,font=\scriptsize,midway] {2};
            \draw[->, shorten >=1pt] (l) -- (ll) node[draw=none,fill=white,font=\scriptsize,midway] {1};
            \draw[->, shorten >=1pt] (rr) -- (rrr) node[draw=none,fill=white,font=\scriptsize,midway] {2};
            \draw[->, shorten >=1pt] (rr) -- (rrl) node[draw=none,fill=white,font=\scriptsize,midway] {1};
            \draw[->, shorten >=1pt] (rl) -- (rlr) node[draw=none,fill=white,font=\scriptsize,midway] {2};
            \draw[->, shorten >=1pt] (rl) -- (rll) node[draw=none,fill=white,font=\scriptsize,midway] {1};
            \draw[->, shorten >=1pt] (lr) -- (lrr) node[draw=none,fill=white,font=\scriptsize, midway] {2};
            \draw[->, shorten >=1pt] (lr) -- (lrl) node[draw=none,fill=white,font=\scriptsize, midway]  {1};
        \end{tikzpicture}
        }
    \end{center}
    \caption{Example optimal decision tree with binary splits, where each node is a state with an action prescription.}
    \label{fig:illexample_dt_trunc}
\end{figure}
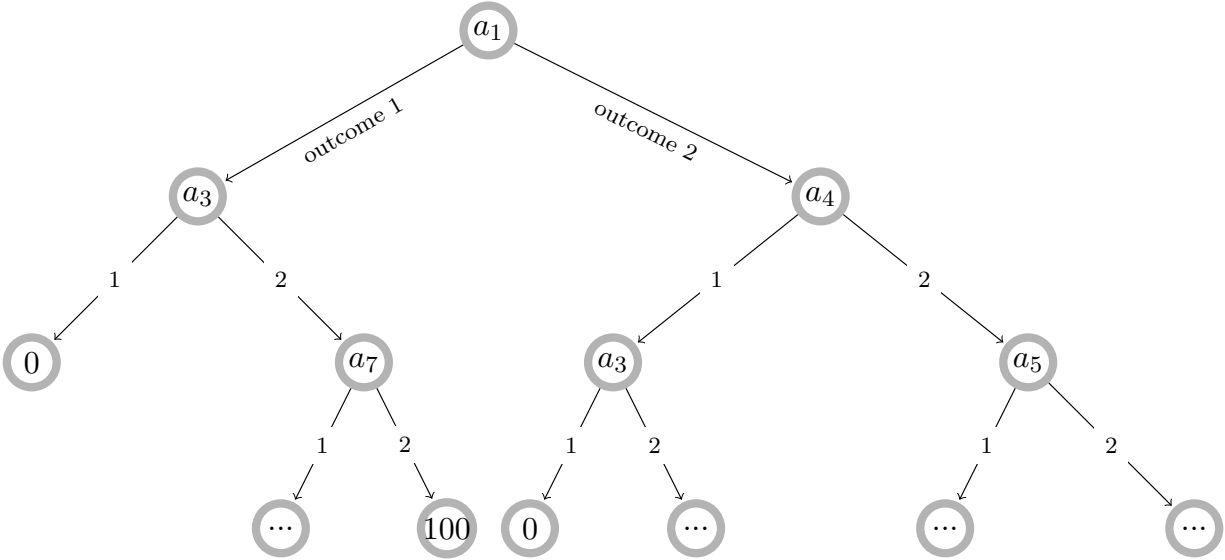

While decision trees are a well-studied model in machine learning, their purpose in that setting is prediction, and each split in the predictive \gls{dt} is a combination of known data features for which the outcome of interest is unknown. In the setting of \gls{coa} generation, the application of \gls{dt}s is for prescription, allowing for the definition of optimal action sequences based on the uncertain outcomes of previous actions. As shown in Figure \ref{fig:illexample_dt_trunc}, each node of a \gls{coa} tree represents a state of the system and prescribes an action, whose outcomes cause the state to change. The state change is represented as a probabilistic move to a child node in the \gls{dt} through one of the outcomes of the action, where another action may be taken. The tree terminates in leaf nodes, where either the action budget is exhausted, the associated states allow for no additional actions, or the agent has accomplished some objective to gain a reward. 

The main contributions of this work are algorithms and optimization formulations to generate globally optimal \gls{dt}s. These \gls{dt}s maximize the expected reward of the strategy, while taking into account the complex interdependencies between actions and outcomes. Using a combination of \gls{dp} and \gls{mio}, our methods scale to problems with large state spaces by using problem-specific information to prune away large subsets of the state space that do not yield progress towards rewards. While we are not the first to suggest reducing the state space via partial action pruning, e.g. in the work of \cite{pinto2014}, we do so without sacrificing global optimality in the final decision model. We use an illustrative example to demonstrate the utility of these methods, and also demonstrate their computational efficiency on a number of randomly generated test cases. Most notably, the presented approach is able to find optimal \gls{dt}s in linear time with respect to the number states explored. 

The framework proposed can generate optimal \gls{dt}s for any problem with the following features:
\begin{enumerate}
    \item An agent takes actions with discrete probabilistic outcomes that change the state of its environment. 
    \item The state captures all relevant information about the environment, as well as information about actions previously taken by the agent.
    \item The decision space is finite, and terminates either when the agent's goals are achieved, or when no additional actions are available.
    \item Actions may have complex interdependencies such as:
    \begin{itemize}
        \item Prerequisites, i.e. actions that must be attempted and result in a specific outcome before another action may be attempted, and
        \item Preclusions, i.e. actions that if attempted and resulting in a specific outcome prevent another action from being attempted. 
    \end{itemize}
\end{enumerate}

The interdependencies between actions can be expressed as a set of logical relationships between actions and associated outcomes, as will be shown mathematically and graphically using an illustrative example in Section~\ref{sec:example}.


\section{Literature Review}
\label{sec:litreview}

In this section, we contrast the use of \gls{dt}s in prediction and \gls{coa} generation settings. We discuss why \gls{dt}s may be preferable to other popular methods such as Markov decision processes and \gls{rl} agents in addressing specific types of decision problems, before discussing the domains where these types of decision problems appear. 

\subsection{Decision trees}
\label{sec:dts}

A \gls{dt} is most commonly encountered in the context of predictive machine learning, in both classification and regression settings\footnote{This section paraphrases author's own work in~\cite{OzturkThesis2022}.}. It is important to differentiate the role of \gls{dt}s in the predictive machine learning setting versus the \gls{coa} generation setting. As a classifier, a \gls{dt} predicts the label of a new data point with known features, or independent variables, but unknown class, or dependent variable, based on \emph{hierarchical splits} on its data features.

In the training phase of a \gls{dt} classifier, the structure of the \gls{dt}, the decision rules at each node of the tree and the labels associated with each leaf are defined. More specifically, starting from data $(\mathbf{X}, \mathbf{y})$, where $\mathbf{X}_i$ is a row of data features associated with class outcome $y_i$, a \gls{dt} $T$ among possible trees $\mathbb{T}$ is chosen by solving the optimization problem
\begin{equation*}
    \underset{T \in \mathbb{T}}{\mathrm{min}}~\mathrm{error}(T, \mathbf{X}, \by) + c_p \cdot \mathrm{complexity}(T).
\end{equation*}
This objective function trades off misclassification error of predicting $\mathbf{y}$ from $\mathbf{X}$ using tree $T$, against tree complexity, which can be measured in terms the depth, breadth and split complexity of the tree, weighted by parameter $c_p$. In general, finding the optimal \gls{dt} via the above formulation is difficult, so approximate methods are common in the academic literature. It is also good practice to perform \emph{cross-validation}, which at its simplest involves splitting available data into training and test sets to avoid overfitting the tree to a specific subset of data.

\gls{dt}s were introduced as a predictive model by~\citet{Breiman1984}, in the form of Classification and Regression Trees (CART). Predictive \gls{dt}s have since proliferated in various forms that improve on a combination of interpretability and predictive power. Some variations are random forests~\citep{breiman2001random}, gradient boosted trees~\citep{chen2016xgboost} and optimal classification trees~\citep{Bertsimas2017}.

The primary similarity between \gls{dt}s in the prediction and \gls{coa} generation settings is in their hierarchical structure. Otherwise, the trees are functionally very different. Firstly, the decisions at each node of the \gls{coa} tree are actions that act on and change the state of the environment; in the prediction setting, the data features are constant and known a priori, and the tree is traversed based on the values of the features without modifying them. The splits in the \gls{coa} tree are probabilistic and based on the outcomes of actions, whereas the splits in the predictive tree are deterministic based on the values of the features. The label of each leaf of a predictive tree is based on the dominant dependent variable value falling in that leaf in the training phase. In the \gls{coa} setting, it is equal to the reward function evaluated at that state, which is dependent on the actions taken in the sequence of nodes in the tree. The predictive \gls{dt} is designed to give the most accurate predictions of a dependent variable through the minimization of misclassification error, while also considering the complexity of the tree. The objective function in the \gls{coa} setting is the expected reward of the leaves, and no regularization or complexity reduction is considered. 

\subsection{Sequential decision making under uncertainty}
\label{sec:decision_literature}

Two popular approaches for sequential decision making under uncertainty are Markov decision processes and reinforcement learning agents. While our work on optimal \gls{dt}s learns and borrows nomenclature from these fields of study, we believe it improves on the limitations of these approaches in addressing certain classes of problems. 

\vspace{5pt} \noindent \textbf{Markov Decision Processes}: The \gls{mdp} is an optimal control problem that involves making sequential decisions to guide a system with stochastic action outcomes towards states that maximize reward. Literature on \gls{mdp}s dates back to the 1950's, when \cite{Bellman1958} applied \gls{dp} to control processes with random noise. The decision problems addressed by \gls{mdp}s have many of the similar features as the ones addressed by \gls{dt}s in this paper, but there are some key distinctions.

The most important distinction is in the concept of state. In both decision models, the state is \emph{memoryless}, meaning that it captures all relevant information about the system and its environment. However, the problems that we attempt to address are ones where the state not only encodes the environment, but also the outcomes of previous actions of the agent. This is critical to our problem definition, since we expect the both the environment and the availability of future actions to evolve in response to the agent's actions. These factors make it impossible for an agent to take actions to revert to a previously visited state. This is in contrast with most \gls{mdp}s, where states describe a system with no memory of past actions, thus allowing for closed state spaces where states can be revisited. Since \gls{mdp} states do not have memory of previous actions by construction, they can only consider simple interdependencies between actions. 

As an extension of the state distinction, a \gls{dt} has advantages in decision problems where there are distinct termination criteria, when no additional actions may be taken or some goals are achieved. While \gls{mdp}s allow for absorbing states~\citep{ermon2014} and action budget constraints~\citep{caramanis2014}, these are considered a special case rather than the norm, and can add significantly to computational complexity. Tree-based search methods have been used to address both MDPs \citep{kocsis2006}, and other games and bandit problems \citep{browne2012}. In these settings, trees were used as a heuristic for approximating one- or multi-step decisions at each state, rather than the fundamental structure for optimizing the entire decision process. 

\vspace{5pt} \noindent \textbf{Reinforcement Learning}: \gls{rl} has become a popular tool for designing policies in environments where maximal performance is valuable, exploration is cheap, and training time is short compared to the time where the learned policy is used \citep{watkins1989}. Deep \gls{rl} has demonstrated success in fields such as computer vision~\citep{le2022}, robotics~\citep{tang2024}, and gaming~\citep{mnih2013}, and is a good candidate for many stochastic decision problems. Here, we review some \gls{rl}-based approaches to \gls{coa} generation and contrast them with a \gls{dt}-based approach.

In a close resemblance to \gls{mdp} approaches, \cite{watkins1989} proposed methods for learning optimal policies in the presence of delayed and uncertain rewards, better known as a Q-learning approach. In Q-learning, the optimal action at any state is the one given by the solution to $\underset{\mathbf{a}}{\rm{argmax}}~Q(\bs, a)$, where $Q: (\bs, a) \xrightarrow{} \mathbb{R}$ is the function that estimates the discounted future reward of taking action $a$ at state $\bs$. In an interesting observation, Watkins does propose the idea that one may follow the courses of action present in a tree of finite depth, but chooses not to pursue due to the intractability of storing a large state space. \cite{barto1983} used a learning system with a search element and a critic element, which has developed into the literature on actor-critic methods. These methods use two agents where the actor agent is in charge of exploring the state space, while the critic element is tasked with determining the value of different actions at different states, much like the Q-function. There is also the more general class of policy search methods~\citep{arul2017}, where the output of the \gls{rl} agent is the likelihood that any action is the best at the current state.

There are substantial challenges in implementing \gls{rl} in the setting defined in Section~\ref{sec:introduction}. The most obvious is the difficulty in strictly enforcing conditional interdependencies between different actions in a setting where the agent learns through experience. These constraints can be satisfied by an \gls{rl} agent either by penalizing the presence of illegal actions in the training procedure, or through post-processing, i.e. by checking whether the actions selected by the optimal policy are valid and choosing the best alternatives otherwise. In most settings however, satisfaction of constraints is still probabilistic \citep{aksaray2021}, which is undesirable in many real-world scenarios.

Another challenge is the need to capture distant rewards with multiple paths to success in the optimal strategy. In the classical cart-pole problem described by ~\cite{barto1983} which has been a staple of the \gls{rl} literature, the reward landscape is a relatively simple function of the state space. However, the problems that we consider have a highly complex relationship between state and reward. It is not obvious that one can train a \gls{rl} agent to successfully exploit disparate pathways to obtain rewards in presence of highly stochastic state transitions, especially if these rewards are temporally distant and in low probability states. This problem has been identified in many papers, but \cite{fan2023} describes it well: ``Learning from sparse rewards is extremely difficult for model-free \gls{rl} algorithms, especially those without intrinsic rewards that struggle to learn from weak gradient signals".

Most importantly, \gls{rl} finds approximate solutions to the \gls{coa} generation problem, and its convergence to a high quality solution is not guaranteed. \cite{sprague2015} demonstrates that deep \gls{rl} is ``brittle" and learns high-quality policies ``within a narrow window" of parameter values, otherwise demonstrating poor performance. \cite{lu2021} argues that it is premature to use deep \gls{rl} in the healthcare setting due to ``its high sensitivity with no clear physiological explanation" to ``input features, embedding model architecture, time discretization, reward function, and random seeds". While there is a wealth of literature on domains where \gls{rl} has succeeded, these often omit the practical difficulties and high computational costs in training reinforcement learning agents. \cite{fan2023} highlights the low learning efficiency of state-of-the-art \gls{rl} algorithms, shows that they are ``fragile facing the hard exploration problems", demonstrating challenges in the kinds of problems that we address in this paper.

\subsection{Applications of decision trees}
\label{sec:applications}

While we borrow much of the nomenclature from \gls{mdp} and \gls{rl} literature, optimal \gls{dt}s are better suited to problems with interdependent actions which require longer and interdependent sequences of actions to reach desired outcomes, such as in wargaming, healthcare and cybersecurity. Optimal \gls{dt}s also provide guarantees of optimality and do not require parameter tuning, which are desirable given the high stakes of decision making in these areas. 

Wargaming is the simulation of conflict that allows decisions makers to explore and analyze scenarios, tactics and strategies for offensive and defensive warfare applications, as well as disaster or biological preparedness.  Such exercises can yield important insights into the effectiveness of tactics and desired concepts of operations.  However, wargames are difficult and require substantial investments to design and execute well \citep{MORSWargaming}. Optimal \gls{dt}s have the potential to add value to  the design of such wargames, by either identifying potentially effective strategies to test, but also in performing criticality analysis for the various actions. 

Healthcare is another domain where optimal \gls{dt}s can make a positive impact. Healthcare providers make recommendations on testing and treatment plans based on available data, and patients decide if and when and which procedures to have done based on their guidance. Such sequential decision making requires determining the next best treatment or diagnosis decision, which changes based on the outcome of previous treatments or updates in the patients' state. Such strategies can be generated through optimization under probabilistic uncertainty. For example, \gls{mdp}s have been applied to determining the optimal cancer treatment plan \citep{CancerPlan}, the optimal management strategy for patients with ischemic heart disease \citep{HAUSKRECHT2000221}, the cost-benefit analysis of deciding when to intervene and operate on patients with hereditary spherocytosis \citep{magni2000}, and optimization of the time to get a living-donor liver donation \citep{LiverDonor}. We hope that our work to extends existing approaches to \gls{coa} generation in the healthcare realm, by both capturing with more granularity the complex relationships between different interventions that naturally arise in healthcare, but also being able to incorporate patient-specific outcome probabilities of different interventions based on individual traits. 

In recent years, cyber space has become an increasingly important dimension of warfare, necessitating the study of advanced cyber warfare analytics~\citep{swallow_cyber}. Identifying optimal strategies for attacking a network is critical for defensive cyber warfare analysis - to both anticipate the actions of attackers, but also take strategic actions to reduce the efficacy of their strategies. \cite{KCCG} present a method for determining the optimal sequencing of phases of a killchain (cyber or physical) with available platforms. A traditional strategy featuring a static action sequence lacks realism due to the probabilistic nature of cyber attacks. For example, attack graphs are often used to represent the privileges acquired by the attacker and the required actions to acquire further privileges \citep{LALLIE2020100219}. While the shortest path through the network could maximize the expectation of acquiring a certain privilege, it disregards the back-up strategies available to an attacker. Thus, using an optimal \gls{dt} to analyze a cyber network could yield richer results.

\section{Illustrative example}
\label{sec:example}

The following example will be used throughout the paper to demonstrate our methods. Table~\ref{fig:illexampletable} contains the parameters of the decision problem. We use the ``(action, outcome)" notation when defining prerequisites and preclusions. For example, a prerequisite $(a_3, 2)$ associated with action $a_7$ means that action $a_7$ can only be taken if action $a_3$ resulted in outcome 2. AND and OR logical constraints are applied as well, in contexts were multiple conditions are present. For example, to be able to take action $a_5$, one must satisfy $(a_4, 2)$, and neither of $(a_3, 1)$ or $(a_3, 2)$. For simplicity, the example problem has actions that have binary outcomes, and each action can only be taken once. Outcomes of certain actions provide reward as indicated in Table~\ref{fig:illexampletable}. However, if multiple rewards are available for any given state, we will define the reward function as the maximum of the applicable rewards based on the action outcomes. In other words, if both action $a_5$ and $a_6$ were taken and both resulted in outcome 2, the reward would be 50.

\begin{table}[h!]
\begin{center}
    \resizebox*{\textwidth}{!}{
        \begin{tabular}{ | c | c | c| c | c | c | c |}
            \hline
            Action & Prerequisite & Preclusion & Cost & Outcome &  Probability & Reward\\
            \hline
            \multirow{2}{*}{$a_1$} &\multirow{2}{*}{--} &\multirow{2}{*}{--} & \multirow{2}{*}{1} & 1 & 0.4 & --\\ 
             &  &  & & 2 & 0.6 & --\\
             \hline
             \multirow{2}{*}{$a_2$} & \multirow{2}{*}{--} &\multirow{2}{*}{--} & \multirow{2}{*}{1} & 1 & 0.4 & --\\
             &  &  & & 2 & 0.6 & --\\  
             \hline
             \multirow{2}{*}{$a_3$} & \multirow{2}{*}{--}  &\multirow{2}{*}{--}& \multirow{2}{*}{1} & 1 & 0.7& --\\  
             &  &  & & 2 & 0.3 & --\\
            \hline
            \multirow{2}{*}{$a_4$} & \multirow{2}{*}{($a_1$,2) OR ($a_3$,2)}  &\multirow{2}{*}{--} & \multirow{2}{*}{1} & 1 & 0.7 & --\\  
             &  &  & & 2 & 0.3 & --\\
            \hline
            \multirow{2}{*}{$a_5$} & \multirow{2}{*}{($a_4$,2)} &\multirow{2}{*}{($a_3$,1) OR ($a_3$,2)} & \multirow{2}{*}{1} & 1 & 0.4 & --\\  
             &  & &  & 2 & 0.6 & 50 \\
            \hline
            \multirow{2}{*}{$a_6$} & \multirow{2}{*}{($a_4$,2) AND ($a_2$,2)} & \multirow{2}{*}{--} & \multirow{2}{*}{1} & 1 & 0.6 & --\\  
             &  & &  & 2 & 0.4 & 10\\
             \hline
             \multirow{2}{*}{$a_7$} & \multirow{2}{*}{($a_3$,2)} & \multirow{2}{*}{--} &  \multirow{2}{*}{1} & 1 & 0.9 & --\\  
             &  & &  & 2 & 0.1 & 100\\
            \hline
        \end{tabular}
    }
    \end{center}
    \caption{Illustrative example parameters.}
    \label{fig:illexampletable}
\end{table}

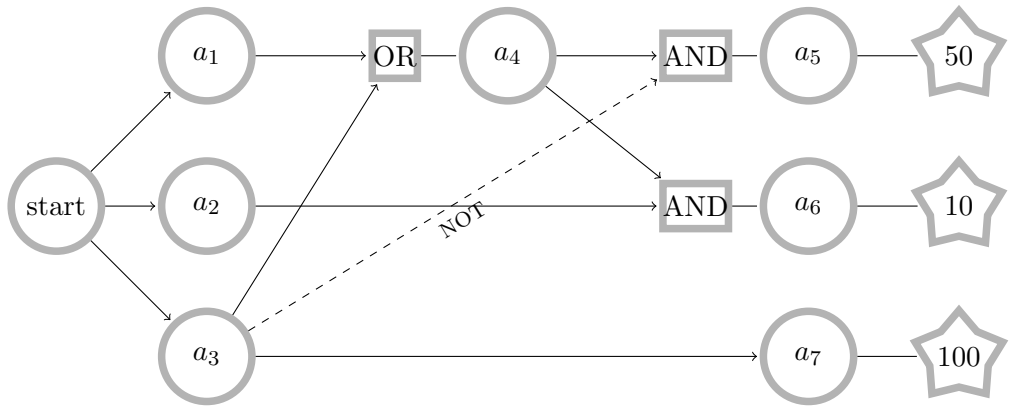
\begin{figure}[h!]
    \begin{center}
        \begin{tikzpicture}
                \node[circle, line width=1mm, inner sep=0pt,
                    minimum size = 12mm,
                    draw=black!30] (start) at (-7,0) {start};
                    
                \node[circle, line width=1mm, inner sep=0pt,
                    minimum size = 12mm,
                    draw=black!30] (1) at (-5,2) {$a_1$};
                    
                \node[circle, line width=1mm, inner sep=0pt,
                minimum size = 12mm,
                draw=black!30] (2) at (-5,0) {$a_2$};
                
                \node[circle, line width=1mm, inner sep=0pt,
                    minimum size = 12mm,
                    draw=black!30] (3) at (-5,-2) {$a_3$};
                    
                \node[rectangle, line width=1mm, inner sep=0pt,
                    minimum size = 6mm,
                    draw=black!30] (OR4) at (-2.5,2) {OR};
                
                \node[circle, line width=1mm, inner sep=0pt,
                    minimum size = 12mm,
                    draw=black!30] (4) at (-1,2) {$a_4$};
                    
                \node[rectangle, line width=1mm, inner sep=0pt,
                    minimum size = 6mm,
                    draw=black!30] (AND5) at (1.5,2) {AND};    
                    
                \node[circle, line width=1mm, inner sep=0pt,
                    minimum size = 12mm,
                    draw=black!30] (5) at (3,2) {$a_5$};
                    
                \node[rectangle, line width=1mm, inner sep=0pt,
                    minimum size = 6mm,
                    draw=black!30] (AND6) at (1.5,0) {AND};
                
                \node[circle, line width=1mm, inner sep=0pt,
                    minimum size = 12mm,
                    draw=black!30] (6) at (3,0) {$a_6$};
                \node[circle, line width=1mm, inner sep=0pt,
                    minimum size = 12mm,
                    draw=black!30] (7) at (3,-2) {$a_7$};
                    
                \node[star, line width=1mm, inner sep =0pt, minimum size = 12mm, draw=black!30] (reward7) at (5,-2) {100};
                
                \node[star, line width=1mm, inner sep =0pt, minimum size = 12mm, draw=black!30] (reward6) at (5,0) {10};
                
                \node[star, line width=1mm, inner sep =0pt, minimum size = 12mm, draw=black!30] (reward5) at (5,2) {50};
                    
                 \draw[->, shorten >=1pt] (start) -- (1); 
                 \draw[->, shorten >=1pt] (start) -- (2);
                 \draw[->, shorten >=1pt] (start) -- (3);
                 \draw[->, shorten >=1pt] (1) -- (OR4);
                 \draw[->, shorten >=1pt] (2) -- (AND6);
                 \draw[->, shorten >=1pt] (3) -- (OR4);
                \draw[dashed, ->, shorten >=1pt] (3) -- (AND5) node[draw=none,fill=none,font=\scriptsize, sloped,midway,below] {NOT};
                \draw[->, shorten >=1pt] (3) -- (7);
                \draw[->, shorten >=1pt] (4) -- (AND5);
                \draw[->, shorten >=1pt] (4) -- (AND6);
                \draw[-, shorten >=1pt] (AND6) -- (6);
                \draw[-, shorten >=1pt] (AND5) -- (5);
                \draw[-, shorten >=1pt] (OR4) -- (4);
                \draw[-, shorten >=1pt] (6) -- (reward6);
                \draw[-, shorten >=1pt] (5) -- (reward5);
                \draw[-, shorten >=1pt] (7) -- (reward7);
        \end{tikzpicture}
    \end{center}
    \caption{Action dependency graph for illustrative example.}
    \label{fig:illexample}
\end{figure}

The prerequisite and preclusion relationships between the actions in Table~\ref{fig:illexampletable} are also shown graphically in Figure \ref{fig:illexample}, as an \gls{adg}. The graph represents prerequisites as solid arrows, where the required outcome of the prerequisite action is outcome 2, and preclusions as dashed arrows. In this example, an action is precluded regardless of the precluding action's outcome. For example, action $a_5$ can only be taken if action $a_4$ resulted in outcome 2 and action $a_3$ has not been taken at all. Finally, the reward gained from an action resulting in outcome 2 is shown in the star-shaped nodes. For example, action $a_5$ with outcome 2 results in a reward of 50. 

\begin{figure}[h!]
    \begin{center}
        \resizebox*{\textwidth}{!}{
        \begin{tikzpicture}
            \node[circle, line width=1mm, inner sep=0pt,
                minimum size = 6mm,
                draw=black!30] (root) at (-1.5,0) {$a_1$};
            \node[circle, line width=1mm, inner sep=0pt,
                minimum size = 6mm,
                draw=black!30] (l) at (-5,-2) {$a_3$};
            \node[circle, line width=1mm, inner sep=0pt,
                minimum size = 6mm,
                draw=black!30] (r) at (2.5,-2) {$a_4$};
            \node[circle, line width=1mm, inner sep=0pt,
                minimum size = 6mm,
                draw=black!30] (ll) at (-7,-4) {0};
            \node[circle, line width=1mm, inner sep=0pt,
                minimum size = 6mm,
                draw=black!30] (lr) at (-3,-4) {$a_7$};
            \node[circle, line width=1mm, inner sep=0pt,
                minimum size = 6mm,
                draw=black!30] (rl) at (0,-4) {$a_3$};
            \node[circle, line width=1mm, inner sep=0pt,
                minimum size = 6mm,
                draw=black!30] (rr) at (5,-4) {$a_5$};
            \node[circle, line width=1mm, inner sep=0pt,
                minimum size = 6mm,
                draw=black!30] (lrr) at (-2,-6) {100};
            \node[circle, line width=1mm, inner sep=0pt,
                minimum size = 6mm,
                draw=black!30] (rll) at (-1,-6) {0};
            \node[circle, line width=1mm, inner sep=0pt,
                minimum size = 6mm,
                draw=black!30] (lrl) at (-4,-6) {$a_2$};
            \node[circle, line width=1mm, inner sep=0pt,
                minimum size = 6mm,
                draw=black!30] (rlr) at (1,-6) {$a_7$};
            \node[circle, line width=1mm, inner sep=0pt,
                minimum size = 6mm,
                draw=black!30] (rrl) at (4,-6) {$a_3$};
            \node[circle, line width=1mm, inner sep=0pt,
                minimum size = 6mm,
                draw=black!30] (rrr) at (7,-6) {$a_3$};
            \node[circle, line width=1mm, inner sep=0pt,
                minimum size = 6mm,
                draw=black!30] (lrlr) at (-3,-8) {$a_4$};
            \node[circle, line width=1mm, inner sep=0pt,
                minimum size = 6mm,
                draw=black!30] (lrll) at (-5,-8) {0};
            \node[circle, line width=1mm, inner sep=0pt,
                minimum size = 6mm,
                draw=black!30] (lrlrr) at (-2,-10) {$a_6$};
            \node[circle, line width=1mm, inner sep=0pt,
                minimum size = 6mm,
                draw=black!30] (lrlrl) at (-4,-10) {0};
            \node[circle, line width=1mm, inner sep=0pt,
                minimum size = 6mm,
                draw=black!30] (lrlrrr) at (-1,-12) {10};
            \node[circle, line width=1mm, inner sep=0pt,
                minimum size = 6mm,
                draw=black!30] (lrlrrl) at (-3,-12) {0};
            \node[circle, line width=1mm, inner sep=0pt,
                minimum size = 6mm,
                draw=black!30] (rlrr) at (2,-8) {100};
            \node[circle, line width=1mm, inner sep=0pt,
                minimum size = 6mm,
                draw=black!30] (rlrl) at (0,-8) {0};
            \node[circle, line width=1mm, inner sep=0pt,
                minimum size = 6mm,
                draw=black!30] (rrll) at (3,-8) {$a_2$};
            \node[circle, line width=1mm, inner sep=0pt,
                minimum size = 6mm,
                draw=black!30] (rrlr) at (5,-8) {$a_7$};
            \node[circle, line width=1mm, inner sep=0pt,
                minimum size = 6mm,
                draw=black!30] (rrrr) at (8,-8) {$a_7$};
            \node[circle, line width=1mm, inner sep=0pt,
                minimum size = 6mm,
                draw=black!30] (rrrl) at (6,-8) {50};
            \node[circle, line width=1mm, inner sep=0pt,
                minimum size = 6mm,
                draw=black!30] (rrlll) at (1,-10) {0};
            \node[circle, line width=1mm, inner sep=0pt,
                minimum size = 6mm,
                draw=black!30] (rrllr) at (3,-10) {$a_6$};
            \node[circle, line width=1mm, inner sep=0pt,
                minimum size = 6mm,
                draw=black!30] (rrlrl) at (4,-10) {0};
            \node[circle, line width=1mm, inner sep=0pt,
                minimum size = 6mm,
                draw=black!30] (rrlrr) at (6,-10) {100};
            \node[circle, line width=1mm, inner sep=0pt,
                minimum size = 6mm,
                draw=black!30] (rrllrl) at (2,-12) {0};
            \node[circle, line width=1mm, inner sep=0pt,
                minimum size = 6mm,
                draw=black!30] (rrllrr) at (4,-12) {10};
            \node[circle, line width=1mm, inner sep=0pt,
                minimum size = 6mm,
                draw=black!30] (rrrrl) at (7,-10) {50};
            \node[circle, line width=1mm, inner sep=0pt,
                minimum size = 6mm,
                draw=black!30] (rrrrr) at (9,-10) {100};

            \draw[->, shorten >=1pt] (root) -- (l) node[draw=none,fill=none,font=\scriptsize, sloped,midway,below] {outcome 1};
            \draw[->, shorten >=1pt] (root) -- (r) node[draw=none,fill=none,font=\scriptsize, sloped,midway,below] {outcome 2};
            \draw[->, shorten >=1pt] (r) -- (rl) node[draw=none,fill=white,font=\scriptsize,midway] {1};
            \draw[->, shorten >=1pt] (r) -- (rr) node[draw=none,fill=white,font=\scriptsize,midway] {2};
            \draw[->, shorten >=1pt] (l) -- (lr) node[draw=none,fill=white,font=\scriptsize,midway] {2};
            \draw[->, shorten >=1pt] (l) -- (ll) node[draw=none,fill=white,font=\scriptsize,midway] {1};
            \draw[->, shorten >=1pt] (rr) -- (rrr) node[draw=none,fill=white,font=\scriptsize,midway] {2};
            \draw[->, shorten >=1pt] (rr) -- (rrl) node[draw=none,fill=white,font=\scriptsize,midway] {1};
            \draw[->, shorten >=1pt] (rl) -- (rlr) node[draw=none,fill=white,font=\scriptsize,midway] {2};
            \draw[->, shorten >=1pt] (rl) -- (rll) node[draw=none,fill=white,font=\scriptsize,midway] {1};
            \draw[->, shorten >=1pt] (lr) -- (lrr) node[draw=none,fill=white,font=\scriptsize, midway] {2};
            \draw[->, shorten >=1pt] (lr) -- (lrl) node[draw=none,fill=white,font=\scriptsize, midway]  {1};
            \draw[->, shorten >=1pt] (lrl) -- (lrlr) node[draw=none,fill=white,font=\scriptsize, midway]  {2};
            \draw[->, shorten >=1pt] (lrl) -- (lrll) node[draw=none,fill=white,font=\scriptsize, midway]  {1};
            \draw[->, shorten >=1pt] (lrlr) -- (lrlrr) node[draw=none,fill=white,font=\scriptsize, midway] {2};
            \draw[->, shorten >=1pt] (lrlr) -- (lrlrl) node[draw=none,fill=white,font=\scriptsize, midway]  {1};
            \draw[->, shorten >=1pt] (lrlrr) -- (lrlrrr) node[draw=none,fill=white,font=\scriptsize, midway]  {2};
            \draw[->, shorten >=1pt] (lrlrr) -- (lrlrrl) node[draw=none,fill=white,font=\scriptsize, midway] {1};
            \draw[->, shorten >=1pt] (rlr) -- (rlrl) node[draw=none,fill=white,font=\scriptsize, midway]  {1};
            \draw[->, shorten >=1pt] (rlr) -- (rlrr) node[draw=none,fill=white,font=\scriptsize, midway] {2};
            \draw[->, shorten >=1pt] (rrr) -- (rrrr) node[draw=none,fill=white,font=\scriptsize, midway] {2};
            \draw[->, shorten >=1pt] (rrr) -- (rrrl) node[draw=none,fill=white,font=\scriptsize, midway]  {1};
            \draw[->, shorten >=1pt] (rrl) -- (rrlr) node[draw=none,fill=white,font=\scriptsize, midway]  {2};
            \draw[->, shorten >=1pt] (rrl) -- (rrll) node[draw=none,fill=white,font=\scriptsize, midway]  {1};
            \draw[->, shorten >=1pt] (rrll) -- (rrlll) node[draw=none,fill=white,font=\scriptsize, midway]  {1};
            \draw[->, shorten >=1pt] (rrll) -- (rrllr) node[draw=none,fill=white,font=\scriptsize, midway] {2};
            \draw[->, shorten >=1pt] (rrlr) -- (rrlrl) node[draw=none,fill=white,font=\scriptsize, midway] {1};
            \draw[->, shorten >=1pt] (rrlr) -- (rrlrr) node[draw=none,fill=white,font=\scriptsize, midway]  {2};
            \draw[->, shorten >=1pt] (rrllr) -- (rrllrl) node[draw=none,fill=white,font=\scriptsize, midway]  {1};
            \draw[->, shorten >=1pt] (rrllr) -- (rrllrr) node[draw=none,fill=white,font=\scriptsize, midway] {2};
            \draw[->, shorten >=1pt] (rrrr) -- (rrrrr) node[draw=none,fill=white,font=\scriptsize, midway] {2};
            \draw[->, shorten >=1pt] (rrrr) -- (rrrrl) node[draw=none,fill=white,font=\scriptsize, midway] {1};
        \end{tikzpicture}
        }
    \end{center}
    \caption{Optimal decision tree for illustrative example, budget = 6.}
    \label{fig:illexample_dt}
\end{figure}
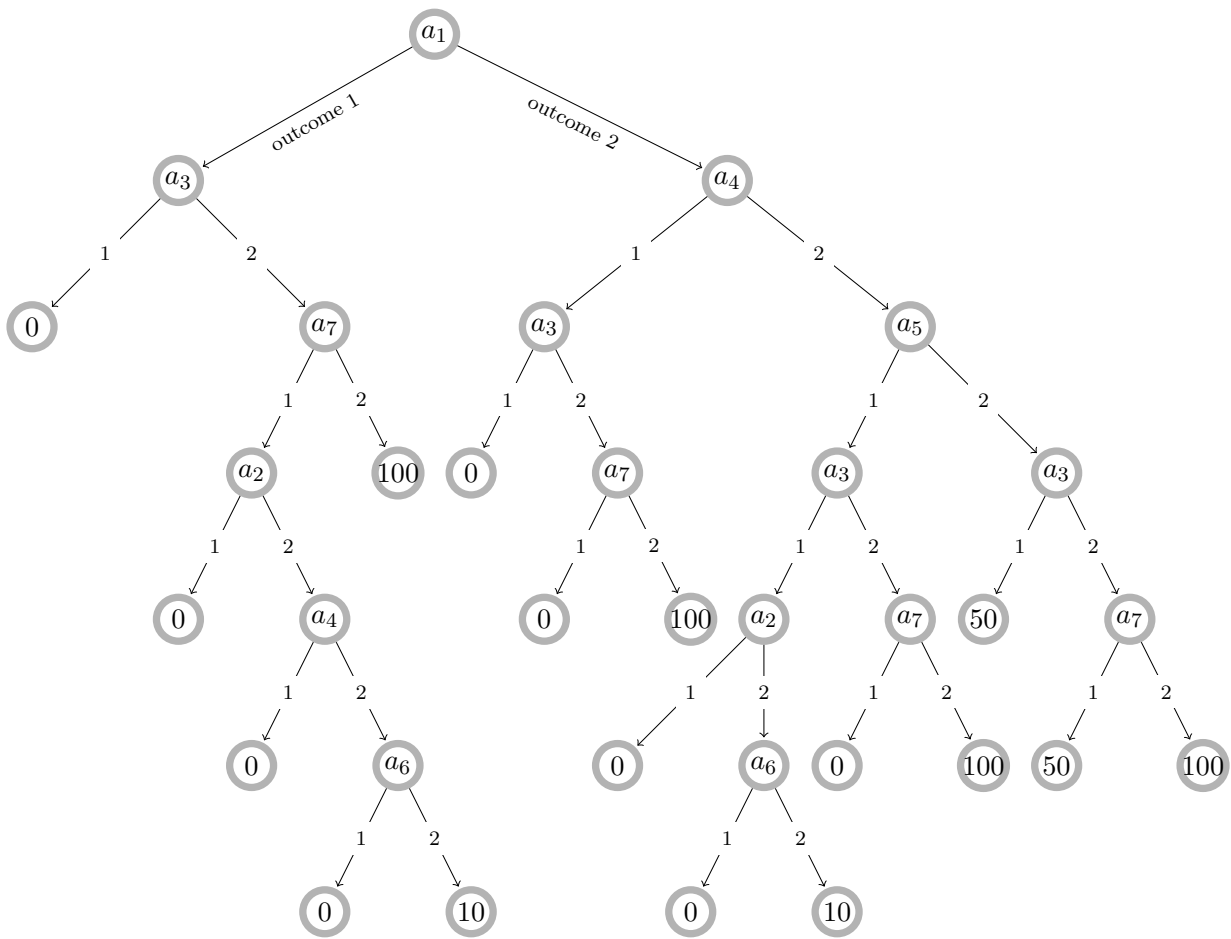

The globally optimal \gls{dt} for this example with a budget of 6 actions is shown in Figure \ref{fig:illexample_dt}. The \gls{dt} shows that action $a_1$ is the optimal action to take first, starting at the root state. If $a_1$ results in outcome 1, action $a_3$ is the next optimal action; otherwise, if $a_1$ results in outcome 2, take action $a_4$. The terminal states indicate the reward received if at that state. For example, if the far right path of exclusively outcome 2's is followed, the reward is 100; alternatively, if the far left path of exclusively outcome 1's is followed, the reward is 0. 

The example highlights the power of the optimal \gls{dt} approach to plan ahead and consider complex tradeoffs between possible strategies. For example, two greedy solutions to this decision problem would have been to choose a strategy with the highest reward ($a_3$ and $a_7$ with a reward of 100 and a probability of $0.03$) or the highest expected reward ($a_1$, $a_4$, and $a_5$ with a reward of 50 and a probability of $0.018$). In presence of action dependencies it is not obvious whether the optimal first action is $a_1$, $a_3$, or even $a_2$, even for a relatively simple decision problem such as this. However, the optimal \gls{dt} captures the benefit of a more certain reward of 50 before making progress towards a reward of 100, due to both the preclusion relationship between $a_3$ and $a_5$ and the disparate rewards of each path. 

\section{Methodology}
\label{sec:methodology}

In order to construct optimal \gls{dt}s, we propose a methodology to explore the state space of the decision problem efficiently, by incorporating problem-specific information about the interdependencies of the actions as well as how rewards are gained. We do this in three stages that we explore in this section, and give an overview of below.
\begin{itemize}
    \item Section~\ref{sec:fullgraph}: We enumerate the subset of all possible states and available actions that have potential for a non-zero reward, called the \emph{full graph}. 
    \item Section~\ref{sec:reducedgraph}: We recursively prune the full graph to determine the \emph{reduced graph}, which is a combination of all subtrees in the full graph that maximize the expectation of reward. 
    \item Section~\ref{sec:optimaltree}: We find the \emph{optimal decision tree} in the reduced graph, breaking ties in the reward expectation via a secondary objective function.
\end{itemize}
While we primarily focus on detailing the process which leverages \gls{dp} and mixed-integer linear optimization, we also include specific software implementation details that can help real-world practitioners execute the theory in this paper. 

\subsection*{Nomenclature}
\label{sec:nomenclature}
In this section, we introduce the mathematical machinery required to generate a \gls{dt}. The variables and functions are defined below:
\begin{itemize}
    \item {State $\bs \in \mathbb{Z}^n$ is a vector representation of the state of the environment. The state uniquely captures all important information about the system of interest.
    \vspace{10pt} \begin{tcolorbox}[title=Illustrative example - State]
        Actions and their outcomes define the entirety of the state space in the example. Each component of the state vector corresponds to an action, and its value denotes whether the action has been taken or its resulting outcome. Specifically, action $a_1$ is represented by its first entry, action $a_2$ is represented by the second, etc.. Each element of the state vector can take one of three values: $0$ (the action has not been taken), $1$ (the action was taken with outcome 1), and $2$ (the action was taken with outcome 2), allowing the state to be represented by combinations of unit vectors. For example, the state $\bs = 2\be_1 +\be_3 + 2 \be_4 = [2,0,1,2,0,0,0]$ represents the state where action $a_1$ and action $a_4$ were both taken with outcome 2, action $a_3$ was taken with outcome 1, and no other actions were taken. 
    \end{tcolorbox}
    }
    \item Root state $\bs_{\rm{root}}$ is the initial state from which \gls{coa} generation starts. 
    \item State domain $\mathbb{S} \subseteq \mathbb{Z}^n$ describes all possible states that can be attained by the system in question, starting from root state $\bs_{\rm{root}}$. 
    \item {Action $r := \{(p_j, \Delta_j), \forall j \in |r|\}$, is the set of probabilities and transition vectors $(p_j, \Delta_j)$ associated with the possible outcomes of the action. Each transition vector $\Delta_j$ is unique for all actions and associated outcomes. We will refer to $(p_j, \Delta_j)$ with or without indices throughout the paper, depending on the context.

    \vspace{10pt} \begin{tcolorbox}[title=Illustrative example - Action]

    Given the state vector definition for the example, the transition vector for action $a_j$ outcome $k$ is $\Delta = k \be_j$, where $\be_j$ is a unit vector with unit $j$'th element. For example, action $a_1$ would be defined as $\{(0.4, \be_1), (0.6, 2\be_1)\}$.
    \end{tcolorbox}
}
    \item {Probability $p_{\bs}$, in a slight abuse of notation, describes probability of reaching state $\bs$, starting from $\bs_{\rm{root}}$.}
    \item {Cost function $c: r \xrightarrow{} \mathbb{R}$ maps each action $r$ to its associated cost.}
    \item {Budget function $B: \bs \xrightarrow{} \mathbb{R}$ returns the remaining budget at state $\bs$. We assume that the budget at the root state ($B(\bs_{\rm{root}})$) is known, and that the budget decreases as more actions are taken.}
    \item {Reward function $\rho: \bs \xrightarrow{} \mathbb{R}$ evaluates the reward achieved at state $\bs$. }
   
    \item {Action function $\mathcal{R}(\bs)$ returns a set of actions $r$ that can be taken at state $\bs$. This function captures whether any actions have conditions that must be met before they can be taken, and whether there is sufficient budget available at $\bs$ to be able to execute the action (i.e. whether $B(\bs) \geq c(r)$).
    \vspace{10pt} \begin{tcolorbox}[title=Illustrative example - Available actions]
    
    Consider the state $\bs = 2\be_1 + \be_4$. The action function $\mathcal{R}(\bs)$ would return $\{a_2, a_3\}$. 
    \end{tcolorbox}
    }
    \item{Decision tree $\mathcal{S}_{\rm{DT}}$ is defined as the set of states contained in a \gls{dt}. }
\end{itemize}

\noindent Starting at a state $\bs$, movement through the state space occurs as follows. Given a choice of action $r \in \mathcal{R}(\bs)$, we obtain the relevant probabilities and state transitions $(p, \Delta) \in r$. The set of children states is generated by adding the state transitions to the current state: $\{\bs + \Delta,~\forall (p, \Delta) \in r\}$. The transition probabilities obey the axioms of probability; given that $\sum_{(p, \Delta) \in r } p = 1$, the probability of reaching state $\bs + \Delta$ is simply the product of the probability of $\bs$ times the action outcome probability: $p_{\bs + \Delta} = p_{\bs} p$. The actions are applied until a leaf node $\bs_{\ell}$ is reached, where no actions are available: $\mathcal{R}(\bs_{\ell}) = \emptyset$. 





Note that we have not discussed how to define $\mathcal{R}(\bs)$; this is problem-specific and would need to be defined for all attainable states. There are no requirements on the mathematical form of $\mathcal{R}(\bs)$ and $\rho(\bs)$. The budget function $B$ must be a decreasing function of the state as additional actions are taken. Additionally, we will assume that the \gls{coa} defined by the \gls{dt} will be executed until a leaf node is reached. Thus, only the leaves of the tree, i.e. each state $\bs$ with $\mathcal{R}(\bs) = \emptyset$, are evaluated when the reward expectation is computed, ensuring that the sum of probabilities of these terminal states will be equal to 1. The resulting reward expectation given the \gls{dt} is 
\begin{equation}
    \mathbb{E}[\rho~|~\mathcal{S}_{\rm{DT}}] = \sum_{\bs \in \mathcal{S}_{\rm{DT}}:~ \mathcal{R}(\bs) = \emptyset} p_{\bs} \rho(\bs).
\end{equation}

We define the following additional machinery, which will become useful as attempt to improve the tractability of generating optimal \gls{dt}s.
\begin{itemize}
   \item{Reward state set $\mathbb{S}_{\rho} := \{\bs:~\bs \in \mathbb{S},~\rho(\bs) > 0\}$: The reward state set is a set that contains all states with a non-zero reward.}
  \item{All action set $\mathcal{R}_{\cup} := \underset{\bs \in \mathbb{S}}{\bigcup} \mathcal{R}(\bs)$: This set contains all possible actions at all possible states.}
  \item{All action-outcome set $\Omega := \{(r, j): r \in \mathcal{R}_{\cup}, (p_j, \Delta_j) \in  r\}$ contains all action-outcome pairs that can change the state of the system.}
  \item{Repeatable action}: A repeatable action is one that can be taken more than once.
\end{itemize}

\subsection{Full graph generation}
\label{sec:fullgraph}

By definition, all possible \gls{dt}s may be enumerated from a state $\bs$ by exhaustively applying each action at each state, one action at a time. However, this is inefficient as the trees in this enumeration will share many states. Instead, if all actions are implemented \emph{in parallel} on all states starting from root state $\bs_{\mathrm{root}}$, and these states are stored, then one can generate a set of all states that can be reached. This we call the \gls{fg}, which is formally the set $\mathcal{S}_{\rm{FG}}$ that is defined by the following recursion: 
\begin{equation}
    \mathcal{S}_{\rm{FG}} = \Big\{\bs: \bs_{\rm{root}} \in \mathcal{S}_{\rm{FG}}; ~\rm{if}~ \bs \in \mathcal{S}_{\rm{FG}},~\bs + \Delta~\in~ \mathcal{S}_{\rm{FG}},~\forall (p, \Delta) \in r,~ \forall r \in \mathcal{R}(\bs) \Big\}
    \label{eq:naiverecursion}
\end{equation}
Algorithm~\ref{alg:Naive_FG} expresses the \gls{fg} recursion in Equation~\eqref{eq:naiverecursion} more explicitly. While we choose to store the \gls{fg} as a set of states, implicitly these states represent a graph through their parent-child relationships. In addition, this set is a graph and not a tree because states may have more than one parent due to different possible orderings of the same set of actions and outcomes.
\begin{algorithm}[h!]
    \caption{Naive Full Graph Generation Algorithm} \label{alg:Naive_FG}
    \begin{algorithmic}[1]
    \State Initialize root state $\bs_{\rm{root}}$.
    \State Initialize full graph state set: $\mathcal{S}_{\rm{FG}} = \{ \bs_{\rm{root}} \}$
    \State{Initialize $\rm{queue} = [\bs_{\rm{root}}]$}.
    \While{$\rm{queue~is~non-empty}$}
        \State $\bs$ = queue.pop().
        \ForAll{actions $r \in \mathcal{R}(\bs)$}
            \ForAll{$(p,\Delta) \in r$}
                \State Generate child $\bs_{\rm{child}} = \bs + \Delta$, with probability $p_{\bs_{\rm{child}}} = p_{\bs}p$.
                \If{$\bs_{\rm{child}} \in \mathcal{S}_{\rm{FG}}$}
                    \State Continue, since child already exists.
                \Else
                    \State Update: $\mathcal{S}_{\rm{FG}} = \mathcal{S}_{\rm{FG}} \cup \{ \bs_{\rm{child}} \}$.
                    \State Add $\bs_{\rm{child}}$ to queue.
                \EndIf
            \EndFor
        \EndFor
    \EndWhile
    \State Return $\mathcal{S}_{\rm{FG}}$. 
    \end{algorithmic}
\end{algorithm}

Given that the set of possible states $\mathbb{S}$ is finite starting from an initial state $\bs_{\rm{root}}$ with a fixed action budget $B(\bs_{\rm{root}})$, it is possible to give a unique numerical index to each state that can be attained in the decision problem by extending the notion of domain to $\mathbb{S}$. This is described in Appendix~\ref{sec:stateindexing} for ease of implementation.

\subsubsection{Accelerating full graph generation}
\label{sec:accel_FG}

The size of the \gls{fg} generated by the naive approach will scale exponentially with the number and repeatability of available actions, and the number of outcomes for each action. Thus the naive approach may be intractable due to the enormity of this state space. In this section, we propose and implement several improvements that reduce the size of the state space without affecting the global optimality of the final \gls{dt}.

\vspace{12pt} \noindent \textbf{Rewarding sets:}
\noindent One approach for greater tractability is to only consider actions at any state that have potential for improved reward. To do so, we need to consider the action-outcome trajectory that can be taken starting from $\bs_{\rm{root}}$ to reach a reward state $\bs_k \in \mathbb{S}_{\rho}$, and how this trajectory evolves as the agent takes actions. Each of these trajectories we call a \emph{rewarding set}, defined mathematically as subset of all action-outcome pairs $\Omega$:
\begin{equation}
    \mathcal{P}(\bs_{\rm{root}}, \bs_{k}) = \Bigg\{\Big\{(r, j) \subseteq \Omega \Big\}: \bs_{k} = \bs_{\rm{root}} + \sum_{(r,j)  \in \mathcal{P}(\bs_{\rm{root}}, \bs_k)} \sum_{(p_j, \Delta_j) \subseteq r} \Delta_j \Bigg\},
    \label{eq:pathrelation}
\end{equation}
where there is a unique set $\mathcal{P}(\bs_{\rm{root}}, \bs_k)$ for each reward state $\bs_k$. As actions are taken that get the agent closer to reaching $\bs_k$ and obtaining reward $\rho(\bs_k)$, the size of the rewarding set leading to $\bs_k$ will decrease, as will be shown shortly.

Note however that, for two rewarding states $\bs_i$ and $\bs_j$, it may be possible to reach $\bs_i$ via a subset of actions and outcomes required to reach $\bs_j$, while simultaneously getting the same or greater reward. In this case, we say that $\bs_i$ \emph{dominates} $\bs_j$. We define the set of dominating reward states $\mathcal{S}_{\rho} \subseteq \mathbb{S}_{\rho}$ as
\begin{equation}
    \mathcal{S}_{\rho}  = \Bigg\{\bs_i: \bs_i \in \mathbb{S}_{\rho},~\nexists s_j \in \mathbb{S}_{\rho}~\rm{s.t.}~\mathcal{P}(\bs_{\rm{root}}, \bs_{j}) \subset \mathcal{P}(\bs_{\rm{root}}, \bs_{i}) \And \rho(\bs_j) \geq \rho(\bs_i) \Bigg\}. 
\end{equation}
While it is possible for the set of all rewarding states $\mathbb{S}_{\rho}$ to be large, the set dominating states $\mathcal{S}_{\rho}$ is usually small. This feature of $\mathcal{S}_{\rho}$ will allow us to parsimoniously keep track of progress towards rewards at each state. From here onward, unless specified otherwise, we will use reward state to refer to a dominating reward state, since these are the most relevant to the algorithm. 

\vspace{10pt} \begin{tcolorbox}[breakable,title=Illustrative example - Rewarding sets]
Recall that various rewards are achieved if action $a_5$, $a_6$, or $a_7$ result in outcome 2. Therefore, the complete reward state set is defined as $\mathbb{S}_{\rho} = \{\bs \in \mathbb{S}: s_5 = 2~|~s_6 = 2~|~s_7 = 2\}$. 
While the non-dominated reward states $\mathcal{S}_{\rho}$ are hard to solve for directly, the rewarding sets $\mathcal{P}(\bs_{\rm{root}}, \bs_k)$ can be identified, one $\bs_k$ at a time, by solving a shortest path algorithm to find the trajectories from the start to the reward nodes, and then excluding the set of paths already considered until all options are exhausted. This process is generalized for any \gls{adg} with and- and or-dependencies as described in Appendix~\ref{sec:rewardingsetsgraph}. For the example problem, the rewarding sets are
\begin{equation}
    \mathbb{P}(\bs_{\rm{root}}) = 
    \begin{cases}
        \{(a_1,2),(a_4,2),(a_5,2)\},\\
        \{(a_1,2),(a_2,2),(a_4,2),(a_6,2)\},\\
        \{(a_2,2),(a_3,2),(a_4,2),(a_6,2)\},\\
        \{(a_3,2),(a_7,2)\},
    \end{cases}
\end{equation}
and can be identified either algorithmically or by observation.
\end{tcolorbox}

\vspace{12pt} \noindent \textbf{Reward-based pruning:}
\noindent Given the above definition of a reward set $\mathcal{P}(\bs_{\rm{root}}, \bs_k)$ from the root state to a dominating reward state $\bs_k \in \mathcal{S}_{\rho}$, at any descendant state $\bs$ that is an arbitrary degree removed from $\bs_{\rm{root}}$, we can store the available rewarding sets based on actions previously taken and the outcomes: $\mathbb{P}(\bs) = \cup_{\bs_k \subseteq \mathcal{S}_{\rho}} \{\mathcal{P}(\bs, \bs_{k})\}$. Given that certain action outcomes may eliminate the viability of a rewarding set, the set of reward states available at $\bs$ are a subset of $\mathcal{S}_{\rho}$. Thus, a child state $\bs_{\rm{child}}$ of a parent state $\bs_{\rm{parent}}$ will inherit a subset of the reward sets available to the parent $\mathbb{P}(\bs_{\rm{parent}})$ using Algorithm \ref{alg:computerewardingsets}.

\begin{algorithm}[h!]
    \caption{Rewarding set computation} \label{alg:computerewardingsets}
    \begin{algorithmic}[1]
        \State {Start at edge $(\bs_{\rm{parent}}, \bs_{\rm{child}})$, where $(r,j)$ is the action-outcome pair that leads from $\bs_{\rm{parent}}$ to  $\bs_{\rm{child}}$.}
        \State Initialize empty set of rewarding sets at the child:  $\mathbb{P}(\bs_{\rm{child}}) = \emptyset$.
        \ForAll{rewarding sets $\mathcal{P}(\bs_{\rm{parent}}, \bs_k) \in \mathbb{P}(\bs_{\rm{parent}})$}
            \If{reward may be increased using the rewarding set: $\rho(\bs) < \rho(\bs_{k})$}
                \If{the action is in the rewarding set: $(r, \_) \in \mathcal{P}(\bs_{\rm{parent}}, \bs_k)$}
                    \If{the action-outcome pair are in rewarding set: $(r, j) \in \mathcal{P}(\bs_{\rm{parent}}, \bs_k)$}
                        \State{Add adjusted rewarding set: $\mathbb{P}(\bs_{\rm{child}}) = \mathbb{P}(\bs_{\rm{child}}) \cup \{\mathcal{P}(\bs_{\rm{parent}}, \bs_k) / (r,j) \}$.}
                    \Else
                        \If{repeatable action and enough budget: \\ $\quad \quad \quad \quad \quad \quad \quad  \quad B(\bs_{\rm{child}}) \geq \sum_{(r', j') \in\mathcal{P}(\bs_{\rm{parent}}, \bs_k)} c(r')$}
                            \State{Add rewarding set: $\mathbb{P}(\bs_{\rm{child}}) = \mathbb{P}(\bs_{\rm{child}}) \cup \{ \mathcal{P}(\bs_{\rm{parent}}, \bs_k) \}.$}
                        \EndIf                    
                    \EndIf
                \Else
                    \If{enough budget: $B(\bs_{\rm{child}}) \geq \sum_{(r', j') \in\mathcal{P}(\bs_{\rm{parent}}, \bs_k)} c(r')$}
                        \State{Add rewarding set: $\mathbb{P}(\bs_{\rm{child}}) = \mathbb{P}(\bs_{\rm{child}}) \cup \{\mathcal{P}(\bs_{\rm{parent}}, \bs_k)\}.$}
                    \EndIf
                \EndIf
            \EndIf
        \EndFor
        \State{Return $\mathbb{P}(\bs_{\rm{child}}).$}
    \end{algorithmic}
\end{algorithm}

In Algorithm~\ref{alg:computerewardingsets}, an underscore outcome in an action outcome-pair denotes that any outcome value is valid for that action, as long as it has been attempted. In addition, the concept of repeatable actions appears; if an action is not repeatable, any outcome of the action that is not the one in the rewarding set is considered to render the rewarding state unreachable. Once $\mathbb{P}(\bs_{\rm{parent}})$ is pruned down into $\mathbb{P}(\bs_{\rm{child}})$, we can simply prune away some of the actions in $\mathcal{R}(\bs_{\rm{child}})$ that do not contribute to reaching a reward state. These rewarding actions are stored in set $\underline{\mathcal{R}}(\bs)$, defined as
\begin{equation}
    \underline{\mathcal{R}}(\bs) = \mathcal{R}(\bs) \bigcap \Bigg\{\underset{~\mathcal{P}(\bs, \bs_k) \in \mathbb{P}(\bs): ~\rho(\bs_k) > \rho(\bs)}{\bigcup} \Big\{ r,~\forall~(r,j) \in \mathcal{P}(\bs, \bs_k) \Big\} \Bigg\}.
    \label{eq:rewardbasedpruning}
\end{equation}
Equation~\eqref{eq:rewardbasedpruning} simply considers all available rewarding sets from $\bs$, and removes the actions in $\mathcal{R}(\bs)$ by checking both whether a given rewarding set has potential for a greater reward than already achieved at $\bs$, but also by checking if an action $r \in \mathcal{R}(\bs)$ has an outcome that is a member of that rewarding set. 

This reward-based pruning approach has several benefits. It can dramatically reduce the number of available actions at a given state, entirely eliminating the exploration of parts of the state space that do not contribute to additional reward. In addition, it allows a user to approximately solve large and perhaps intractable decision problems. Since one may add rewarding sets $\mathcal{P}(\bs_{\rm{root}}, \bs_k)$ incrementally to the decision problem by only considering a subset of the rewarding states $\bs_k \in \mathcal{S}_{\rho}$ at a time, one may generate a sub-optimal \gls{coa} using a subset of the \gls{fg}, thus find a lower bound on the optimal \gls{coa} including all rewarding states. While we solve to global optimality and thus do not pursue this approach in this paper, we believe it is an important contribution that allows practitioners to get near-optimal solutions to otherwise intractable decision problems without surpassing their computational budget. 

It is important to reiterate that the reason for reward-based pruning is to save time and memory by reducing state space enumeration. In this context, the need to store $\mathbb{P}(\bs)$ will increase time and memory requirements, potentially offsetting some reduction in time and memory savings due to fewer states explored. However, this is not a concern with a clever implementation of the algorithm. Once the children of a state $\bs$ are initialized, the rewarding sets $\mathbb{P}(\bs)$ may be forgotten in order to reduce the memory requirements. In addition, a depth-first implementation of the approach, as will be showcased in Algorithm~\ref{alg:fullgraph}, further minimizes the memory requirements by storing as few unexplored states in the graph as possible, thus further minimizing stored rewarding set information. Even without such considerations, the reduction of actions from $\mathcal{R}(\bs)$ to $\underline{\mathcal{R}}(\bs)$ via Equation~\eqref{eq:rewardbasedpruning} allows for substantial speed-ups beyond the time penalty required to compute $\mathcal{P}(\bs)$, except in the most adversarial cases.

\vspace{10pt} \begin{tcolorbox}[breakable,title=Illustrative example - Pruning]

\begin{center}
    \resizebox*{0.7\textwidth}{!}{
        \begin{tikzpicture}
                \node[ellipse, line width=1mm, inner sep=0pt,
                    minimum height = 12mm, minimum width = 24mm,
                    draw=black!30] (root) at (-3, 0) {$\Vec{0}$};
                    
                \node[ellipse, line width=1mm, inner sep=0pt,
                    minimum height = 12mm, minimum width = 24mm,
                    draw=black!30] (l_1) at (-6, -2) {$\be_1$};
                
                \node[ellipse, line width=1mm, inner sep=0pt,
                    minimum height = 12mm, minimum width = 24mm,
                    draw=black!30] (c_1) at (-3, -2) {$2\be_1$};
                
                \node[ellipse, line width=1mm, inner sep=0pt,
                    minimum height = 12mm, minimum width = 24mm,
                    draw=black!30] (r_1) at (0, -2) {...};
                    
                \node[ellipse, line width=1mm, inner sep=0pt,
                    minimum height = 12mm, minimum width = 24mm,
                    draw=black!30] (l_2) at (-6, -4) {...};
                
                \node[ellipse, line width=1mm, inner sep=0pt,
                    minimum height = 12mm, minimum width = 24mm,
                    draw=black!30] (c_2) at (-3, -4) {$2\be_1 + \be_4$};
                
                \node[ellipse, line width=1mm, inner sep=0pt,
                    minimum height = 12mm, minimum width = 24mm,
                    draw=black!30] (r_2) at (0, -4) {$2\be_1 + 2\be_4$};
                    
                \node[ellipse, line width=1mm, inner sep=0pt,
                    minimum height = 12mm, minimum width = 24mm,
                    draw=black!30] (rr_2) at (3, -4) {...};
            
                \draw[->, shorten >=1pt] (root) -- (l_1);
                \draw[->, shorten >=1pt] (root) -- (c_1);
                \draw[->, shorten >=1pt, dashed] (root) -- (r_1);
                
                \draw[->, shorten >=1pt, dashed] (l_1) -- (l_2);
                \draw[->, shorten >=1pt] (c_1) -- (c_2);
                \draw[->, shorten >=1pt] (c_1) -- (r_2);
                \draw[->, shorten >=1pt, dashed] (c_1) -- (rr_2);
        \end{tikzpicture}
    }
\end{center}

Let us apply the reward-based pruning procedure to determine the available actions at the states in the graph above. The following table summarizes $\mathcal{R}(\bs)$, $\mathbb{P}(\bs)$ and the resulting pruned $\underline{\mathcal{R}}(\bs)$. Note that the rewarding sets at a given state are the actions that have yet to be taken to reach reward and thus get smaller as more actions are taken. \\

\begin{center}
\begin{tabular}{|c|c c c|}
\hline
     State & $\mathcal{R}(\bs)$ & $\mathbb{P}(\bs)$ & $\underline{\mathcal{R}}(\bs)$ \\
     \hline
     \multirow{4}{*}{$\Vec{0}$} & \multirow{4}{*}{$\{a_1, a_2, a_3 \}$} & $\{\{(a_1,2),(a_4,2),(a_5,2)\},$ & \multirow{4}{*}{$\{a_1, a_2, a_3 \}$}\\
     & & $\{(a_1,2),(a_2,2),(a_4,2),(a_6,2)\},$ & \\
     & & $\{(a_2,2),(a_3,2),(a_4,2),(a_6,2)\},$ & \\
     & & $\{(a_3,2),(a_7,2)\}\}$ & \\
     \hline
     \multirow{2}{*}{$\be_1$}& \multirow{2}{*}{$\{a_2, a_3\}$} & $\{(a_2,2),(a_3,2),(a_4,2),(a_6,2)\},$ & \multirow{2}{*}{$\{a_2, a_3\}$}\\
     & & $\{(a_3,2),(a_7,2)\}\}$ &\\
     \hline
     \multirow{4}{*}{$2\be_1$}& \multirow{4}{*}{$\{a_2, a_3, a_4 \}$} & $\{\{(a_4,2), (a_5,2)\},$ & \multirow{4}{*}{$\{a_2, a_3, a_4\}$}\\
     & & $\{(a_2,2),(a_4,2),(a_6,2)\},$ & \\
     & & $\{(a_2,2),(a_3,2),(a_4,2),(a_6,2)\},$ &\\
     & & $\{(a_3,2),(a_7,2)\}\}$ & \\
     \hline
     $2\be_1 + \be_4$ & $\{a_2, a_3\}$ & $\{\{(a_3,2),(a_7,2)\}\}$ & $\{a_3\}$\\
     \hline
     \multirow{4}{*}{$2\be_1 + 2\be_4$}& \multirow{4}{*}{$\{a_2, a_3, a_5\}$} & $\{\{(a_5,2)\},$ & \multirow{4}{*}{$\{a_2, a_3, a_5\}$}\\
     & & $\{(a_2,2),(a_6,2)\},$ & \\
     & & $\{(a_2,2),(a_3,2),(a_6,2)\},$ & \\
     & & $\{(a_3,2),(a_7,2)\}\}$ & \\
     \hline
\end{tabular}
\end{center}
At the root state $\Vec{0}$ all rewarding sets are feasible, meaning that $\mathcal{R}(\Vec{0}) = \underline{\mathcal{R}}(\Vec{0})$. Similarly, the state $2\be_1$ inherits all rewarding sets from its parent, the root, and therefore its pruned actions are the same as its initial set of actions. 

More notably, having taken action $a_1$ with outcome 1, the state $\be_1$ does not inherit the two rewarding sets that require action $a_1$ with outcome 2 and are thus no longer available. However, this does not prune the available actions because both action $a_2$ and action $a_3$ are present in the remaining rewarding sets. The state $2\be_1 + \be_4$ inherits only one rewarding set from its parent $2\be_1$ because the action-outcome $(a_4, 2)$ was present in the other three rewarding sets. Because action $a_2$ is not present in the remaining rewarding set, we can prune it away, saving computational effort. 

\end{tcolorbox}

\vspace{12pt} \noindent \textbf{Accelerated full graph algorithm:} Algorithm~\ref{alg:fullgraph} combines the naive implementation of \gls{fg} generation (Algorithm~\ref{alg:Naive_FG}) with rewarding set computations (Algorithm~\ref{alg:computerewardingsets}) and the associated determination of the effective actions  (Equation~\ref{eq:rewardbasedpruning}). This way, it reduces the size of the \gls{fg} with zero impact on the score of the optimal tree. 

\begin{algorithm}[h!]
    \caption{Accelerated Full Graph Generation Algorithm} \label{alg:fullgraph}
    \begin{algorithmic}[1]
    \State Initialize root state $\bs_{\rm{root}}$.
    \State Determine actions at root state: $\mathcal{R}(\bs_{\rm{root}})$. 
    \State {Compute rewarding sets: $\mathbb{P}(\bs_{\rm{root}}) = \{\mathcal{P}(\bs_{\rm{root}}, \bs_k), \forall \bs_k \in \mathcal{S}_{\rho}\}$.}
    \State {Prune actions (Eq.~\ref{eq:rewardbasedpruning}): $\mathcal{R}(\bs_{\rm{root}}) \xrightarrow{} \underline{\mathcal{R}}(\bs_{\rm{root}})$}.
    \State Initialize full graph set: $\mathcal{S}_{\rm{FG}} = \{\bs_{\rm{root}}\}$.
    \State{Initialize $\rm{queue} = [\bs_{\rm{root}}]$}.
    \While{$\rm{queue~is~non-empty}$}
        \State $\bs$ = queue.pop().
        \ForAll{actions $r \in \underline{\mathcal{R}}(\bs)$}
            \ForAll{$(p,\Delta) \in r$}
                \State Candidate child $\bs_{\rm{child}} = \bs + \Delta, p_{\bs_{\rm{child}}} = p_{\bs} p$.
                \If{$\bs_{\rm{child}} \in \mathcal{S}_{\rm{FG}}$}
                    \State Continue, since child already exists.
                \Else
                    \State Update: $\mathcal{S}_{\rm{FG}} = \mathcal{S}_{\rm{FG}} \cup \{\bs_{\rm{child}} \}$
                    \State Determine actions at child state: $\mathcal{R}(\bs_{\rm{child}})$. 
                    \State {Compute rewarding sets: $\mathbb{P}(\bs_{\rm{child}}) \subseteq \mathbb{P}(\bs)$.}
                    \State {Prune actions (Eq.~\ref{eq:rewardbasedpruning}): $\mathcal{R}(\bs_{\rm{child}}) \xrightarrow{} \underline{\mathcal{R}}(\bs_{\rm{child}})$}.
                    \State Add $\bs_{\rm{child}}$ to queue.
                \EndIf
            \EndFor
        \EndFor
    \EndWhile
    \State Return $\mathcal{S}_{\rm{FG}}$. 
    \end{algorithmic}
\end{algorithm}

\subsection{Reduced graph generation}
\label{sec:reducedgraph}

Now that we have enumerated the possible strategies that lead to reward states, we devise methods to reduce $\mathbb{S}_{\rm{FG}}$ to the subset of states that exist in the union of possible optimal strategies, i.e. in the union of the optimal \gls{dt}s. We call this subset the \emph{\gls{rg}}, represented by $\mathcal{S}_{\rm{RG}} \subseteq \mathcal{S}_{\rm{FG}}$. To define $\mathcal{S}_{\rm{RG}}$, we need to first compute the maximum reward expectation of each state in the graph. This can be done by recursively comparing the reward expectation of all subtrees starting at each state to each other. This process defines function $\Phi$, which gives the maximum reward expectation possible at any state $\bs$ in the \gls{fg}. The process for defining the map $\Phi: \bs \xrightarrow{} \mathbb{R}$ is shown in the recursion in Equation~\eqref{eq:scorerecursion}.
\begin{equation}
    \Phi(\bs) = 
    \begin{cases}
    \rho(\bs),~\rm{if}~\underline{\mathcal{R}}(\bs) = \emptyset, \\
    \underset{r \in \underline{\mathcal{R}}(\bs)}{\max} \Bigg(\underset{(p, \Delta) \in r}{\sum} p \Phi(\bs + \Delta) \Bigg),~\rm{otherwise}.
    \end{cases}
    \label{eq:scorerecursion}
\end{equation}
$\Phi(\bs)$ is simply equivalent to the reward $\rho(\bs)$ if no further actions are available at $\bs$. Otherwise, $\Phi(\bs)$ is the reward expectation achieved by taking an action $r \in \underline{\mathcal{R}}(\bs)$ that maximizes the expectation of $\Phi$ at the resulting children states. In practice, computing the score of a given state is a \gls{dp}, and once the score is evaluated at $\bs$, it must also be true that the scores have been computed for all children states of $\bs$ of any degree. These can be efficiently stored in memory for future queries without having to recompute the recursion.

The \gls{rg} is simply the \gls{fg} that has had all dominated actions pruned. This is done using the recursion in Algorithm~\ref{alg:reducedgraph}.
\begin{algorithm}
    \caption{Reduced Graph Generation Algorithm}\label{alg:reducedgraph}
    \begin{algorithmic}[1]
   \State Initialize root state $\bs_{\rm{root}}$.
   \State Compute $\Phi(\bs_{\rm{root}})$ via Eq.~\ref{eq:scorerecursion}.
    \State Initialize reduced graph state set: $\mathcal{S}_{\rm{RG}} = \{\bs_{\rm{root}}\}$. 
    \State{Initialize $\rm{queue} = [\bs_{\rm{root}}]$}.
    \While{$\rm{queue~is~non-empty}$}
        \State $\bs$ = queue.pop().
        \State Get best action set, $\mathcal{R}^*(\bs) = \{r: \sum_{(p,\Delta) \in r} p\Phi(\bs + \Delta) = \Phi(\bs), r \in \underline{\mathcal{R}}(\bs)\}$
        \ForAll{actions $r \in \mathcal{R}^*(\bs)$}
            \ForAll{$(p,\Delta) \in r$}
                \State Candidate child $\bs_{\rm{child}} = \bs + \Delta, p_{\bs_{\rm{child}}} = p_{\bs} p$.
                \If{$\bs_{\rm{child}} \in \mathcal{S}_{\rm{RG}}$}
                    \State Continue, since child already exists.
                \Else
                    \State Update: $\mathcal{S}_{\rm{RG}} = \mathcal{S}_{\rm{RG}} \cup \{\bs_{\rm{child}} \}$
                    \State Add $\bs_{\rm{child}}$ to queue.
                \EndIf
            \EndFor
        \EndFor
    \EndWhile
    \end{algorithmic}
\end{algorithm}
Since we had already generated and stored the \gls{fg} and the actions at each state $\underline{\mathcal{R}}(\bs)$, we recall the various available actions starting at the root state, and make decisions among them, only choosing to execute the most effective actions at each state.  As an additional note, in Algorithm~\ref{alg:reducedgraph}, there is no need to use the rewarding sets, since for any state $\bs$, the reduced action set $\underline{\mathcal{R}}(\bs)$ contains only the actions that would lead to increased reward.

The computation of the root score in line 2 of Algorithm~\ref{alg:reducedgraph} implies that the scores of all nodes is computed prior to starting to add states to the \gls{rg}. This score computation is necessary because of the breadth-first, the top-down construction of the \gls{rg}. This design ensures that only the nodes that maximize the score are re-explored, but requires information about expected reward to be communicated from the terminal nodes of the \gls{fg} all the way to the root state. 

\vspace{10pt} \begin{tcolorbox}[breakable,title=Illustrative example - Reduced graph generation]

To compute $\Phi([2,0,2,2,1,0,0])$, we have to consider its subgraph and determine $\Phi(\bs)$ for all states in that subgraph, as shown below. Note that in the subgraph, states may have parents originating from other states that are not shown.

\begin{center}
    \resizebox*{\textwidth}{!}{
        \begin{tikzpicture}
                \node[ellipse, line width=1mm, inner sep=0pt,
                    minimum height = 12mm, minimum width = 24mm,
                    draw=black!30] (root) at (0, 0) {$[2,0,2,2,1,0,0]$};

                \node[ellipse, line width=1mm, inner sep=0pt,
                    minimum height = 12mm, minimum width = 24mm,
                    draw=black!30] (ll_1) at (-6.5, -2) {$[2,1,2,2,1,0,0]$};
                
                \node[ellipse, line width=1mm, inner sep=0pt,
                    minimum height = 12mm, minimum width = 24mm,
                    draw=black!30] (l_1) at (-2.5, -2) {$[2,2,2,2,1,0,0]$};
                
                \node[ellipse, line width=1mm, inner sep=0pt,
                    minimum height = 12mm, minimum width = 24mm,
                    draw=black!30] (r_1) at (2.5, -2) {$[2,0,2,2,1,0,1]$};
                    
                \node[ellipse, line width=1mm, inner sep=0pt,
                    minimum height = 12mm, minimum width = 24mm,
                    draw=black!30] (rr_1) at (6.5, -2) {$[2,0,2,2,1,0,2]$};

                \node[ellipse, line width=1mm, inner sep=0pt,
                    minimum height = 12mm, minimum width = 24mm,
                    draw=black!30] (lll_2) at (-12, -4) {$[2,1,2,2,1,0,1]$};
                
                \node[ellipse, line width=1mm, inner sep=0pt,
                    minimum height = 12mm, minimum width = 24mm,
                    draw=black!30] (ll_2) at (-8, -4) {$[2,1,2,2,1,0,2]$};

                \node[ellipse, line width=1mm, inner sep=0pt,
                    minimum height = 12mm, minimum width = 24mm,
                    draw=black!30] (l_2) at (-3.5, -4) {$[2,2,2,2,1,1,0]$};
                    
                \node[ellipse, line width=1mm, inner sep=0pt,
                    minimum height = 12mm, minimum width = 24mm,
                    draw=black!30] (r_2) at (0.5, -4) {$[2,2,2,2,1,2,0]$};
                    
                \node[ellipse, line width=1mm, inner sep=0pt,
                    minimum height = 12mm, minimum width = 24mm,
                    draw=black!30] (rr_2) at (5, -4) {$[2,2,2,2,1,0,1]$};
                
                \node[ellipse, line width=1mm, inner sep=0pt,
                    minimum height = 12mm, minimum width = 24mm,
                    draw=black!30] (rrr_2) at (9, -4) {$[2,2,2,2,1,0,2]$};

                \draw[->, shorten >=1pt] (root) -- (ll_1) node[draw=none,fill=none,midway,below, shift={(1.5,0.3)}] {$a_2$};
                \draw[->, shorten >=1pt] (root) -- (l_1);
                \draw[->, shorten >=1pt] (root) -- (r_1);
                \draw[->, shorten >=1pt] (root) -- (rr_1) node[draw=none,fill=none,midway,below, shift={(-1.5,0.3)}] {$a_7$};
                
                \draw[->, shorten >=1pt] (ll_1) -- (lll_2) node[draw=none,fill=none,midway,right, shift={(1.0, 0.0)}] {$a_7$};
                \draw[->, shorten >=1pt] (ll_1) -- (ll_2);
                \draw[->, shorten >=1pt] (l_1) -- (l_2) node[draw=none,fill=none,midway,right, shift={(0.7, 0.)}] {$a_6$};
                \draw[->, shorten >=1pt] (l_1) -- (r_2);
                \draw[->, shorten >=1pt] (l_1) -- (rr_2);
                \draw[->, shorten >=1pt] (l_1) -- (rrr_2) node[draw=none,fill=none,midway,below] {$a_7$};
        \end{tikzpicture}
    }
\end{center}
The scores and optimal action sets are summarized in the following table:

\begin{center}
\begin{tabular}{|c|c c c|}
\hline
     State & $\underline{\mathcal{R}}(\bs)$ & Score & $\mathcal{R}^*(\bs)$\\
     \hline
     $[2,0,2,2,1,0,0]$ & $\{a_2, a_7\}$ & $\max(0.1, 0.1) = 0.1$ & $\{a_2, a_7\}$\\
     \hline
     $[2,1,2,2,1,0,0]$ & $\{a_7\}$ & $\max(0.1) = 0.1$ & $\{a_7\}$\\
     $[2,2,2,2,1,0,0]$ & $\{a_6, a_7\}$ & $\max(0.04, 0.1) = 0.1$ & $\{a_7\}$\\
     $[2,0,2,2,1,0,1]$ & $\emptyset$ & 0 & $\emptyset$ \\
     $[2,0,2,2,1,0,2]$ & $\emptyset$ & 1 & $\emptyset$ \\
     \hline
     $[2,1,2,2,1,0,1]$ & $\emptyset$ & 0 & $\emptyset$ \\
     $[2,1,2,2,1,0,2]$ & $\emptyset$ & 1 & $\emptyset$ \\
     $[2,2,2,2,1,1,0]$ & $\emptyset$ & 0 & $\emptyset$ \\
     $[2,2,2,2,1,2,0]$ & $\emptyset$ & 0.1 & $\emptyset$ \\
     $[2,2,2,2,1,0,1]$ & $\emptyset$ & 0 & $\emptyset$ \\
     $[2,2,2,2,1,0,2]$ & $\emptyset$ & 1 & $\emptyset$ \\
     \hline
\end{tabular}
\end{center}

There are multiple optimal actions to take at the state $[2,0,2,2,1,0,0]$, namely $a_2$ and $a_7$, so the optimal subtrees associated with both actions would be present in the reduced graph. However, for state $[2,2,2,2,1,0,0]$, action $a_6$ is not in $\mathcal{R}^*(\bs)$, and the children associated with taking action $a_6$ would not be included in the reduced graph. 
\end{tcolorbox}

\subsubsection{Alternative approach: linear programming}


An alternate method for determining the \gls{rg} is to formulate and solve a \gls{lo} problem. While our experiments have shown that the \gls{dp} approach is computationally more efficient, the \gls{lo} approach enables adversarial extensions of the decision problem. In this robust optimization context, an agent may need to consider a policy that maximizes reward while being robust to limited adversarial perturbations to the outcome probabilities of its actions. In addition, using \gls{lo} could be more efficient for initial development of our \gls{dt} algorithm on new applications. Given that \gls{dp}s can be complicated to implement in many scenarios, the \gls{lo} could be used to validate the \gls{dp}-based \gls{rg} generation. Thus, we give the \gls{lo} formulation of the score recursion in Equation~\eqref{eq:scorerecursion}:
\begin{equation}
    \begin{split}
        \underset{\Phi}{\rm{minimize}} \quad & \sum_{\bs \in \mathcal{S}_{\rm{FG}}}\Phi(\bs) \\ 
    \rm{subject~to~}\quad & \sum_{(p, \Delta) \in r} p\Phi(\bs +\Delta) \leq \Phi(\bs),~\forall r \in \underline{\mathcal{R}}(\bs), \forall \bs \in \mathcal{S}_{\rm{FG}} \setminus \mathcal{S}_L, \\ 
    & \Phi(\bs)=\rho(\bs), \forall \bs \in \mathcal{S}_L, \\
    & \mathcal{S}_L := \{\bs \in \mathcal{S}_{\rm{FG}}: \underline{\mathcal{R}}(\bs) = \emptyset\}, \\
    & \Phi(\bs) \in [0,1], \forall \bs \in \mathcal{S}_{\rm{FG}}.
    \end{split} \label{eq:LOmethod}
\end{equation}

\noindent Given the full graph $\mathcal{S}_{\rm{FG}}$, let $\Phi(\bs)$ be the score variable for each state $\bs \in \mathcal{S}_{\rm{FG}}$. $\mathcal{S}_L$ is the set of states with no children ($\underline{\mathcal{R}}(\bs) = \emptyset$); the score at those states is constrained to be equal to the reward at that state. For all other states, the score is lower-bounded by the expected reward given a possible action, for all possible actions. By minimizing the sum of the scores at all states, the scores will assume the value of their largest lower-bound, thus maximizing the expected reward at the root. At any given state, the constraints which have no slack, corresponding to non-zero dual variables, point to the set of optimal actions at that state. 

This formulation naturally lends itself to an adversarial extension which the authors hope to explore in future work. Consider the case where the \gls{dt} is the optimal strategy for executing a cyber attack against a network. The network's owner might have an arsenal of tools to make the attacker's actions more expensive or more difficult. Thus, the adversarial extension would need to identify which actions to interdict to minimize the attacker's expected reward over all possible strategies in from the \gls{fg}. Because the \gls{lo} is already formulated as a minimization problem, it is amenable to the incorporation of the interdictor's decisions and budget, unlike the \gls{dp} method.

\subsection{Tree selection}
\label{sec:optimaltree}

Tree selection requires a final recursion, given in Algorithm~\ref{alg:optimaltree}, traversing the \gls{rg} so that only one action is taken for all attainable states starting at $\bs_{\rm{root}}$. Given that all actions in the \gls{rg} have the same reward expectation, this can be done based on a secondary objective, which we simply leave as $f: \bs \xrightarrow{} \mathbb{R}$.
\begin{algorithm}[h!]
    \caption{Optimal Tree Generation Algorithm}\label{alg:optimaltree}
    \begin{algorithmic}[1]
   \State Initialize root state $\bs_{\rm{root}}$.
    \State Compute $f(\bs), ~\forall \bs \in \mathcal{S}_{\rm{RG}}$. 
    \State Initialize decision tree state set: $\mathcal{S}_{\rm{DT}} = \{\bs_{\rm{root}}\}$. 
    \State{Initialize $\rm{queue} = [\bs_{\rm{root}}]$}.
    \While{$\rm{queue~is~non-empty}$}
        \State $\bs$ = queue.pop().
        \State Choose among iso-reward actions: $r^* = \underset{r \in \mathcal{R}^*(\bs)}{\rm{argmin}} \sum_{(p,\Delta) \in r} pf(\bs + \Delta)$. 
        \ForAll{$(p,\Delta) \in r^*$}
            \State Child $\bs_{\rm{child}} = \bs + \Delta, p_{\bs_{\rm{child}}} = p_{\bs} p$.
            \State Update: $\mathcal{S}_{\rm{DT}} = \mathcal{S}_{\rm{DT}} \cup \{ \bs_{\rm{child}} \}$
            \State Add $\bs_{\rm{child}}$ to queue.
        \EndFor
    \EndWhile
    \State
    Return $\mathcal{S}_{\rm{DT}}$.
    \end{algorithmic}
\end{algorithm}
Some choices of $f$ may minimize the number of nodes or depth of the optimal tree starting from $\bs$, or bias the optimum towards \gls{dt}s based on orderings of actions. Given that all subtrees in the \gls{rg} have the same score, one may simply pick actions at random as well. The complexity of computing $f(\bs)$ will have an impact on the speed of Algorithm~\ref{alg:optimaltree}. By default, we will assume that $f $ is the function that returns the number of nodes of the optimal tree originating from a state. Mathematically, 
\begin{equation}
    f(\bs) = \begin{cases}
        1,~&\rm{if}~\mathcal{R}^*(\bs) = \emptyset, \\
        \underset{r \in \mathcal{R}^*(\bs)}{\min} \sum_{(p, \Delta) \in r} f(\bs + \Delta),~&\rm{otherwise},
    \end{cases}
\end{equation}
which is a recursion similar to the score recursion in Equation~\eqref{eq:scorerecursion}.

\vspace{10pt} \begin{tcolorbox}[breakable,title=Illustrative example - Optimal decision tree]
Consider the reduced graph originating from state $[2,0,2,2,1,0,0]$ shown below. 

\begin{center}
    \resizebox*{\textwidth}{!}{
        \begin{tikzpicture}
                \node[ellipse, line width=1mm, inner sep=0pt,
                    minimum height = 12mm, minimum width = 24mm,
                    draw=black!30] (root) at (0, 0) {$[2,0,2,2,1,0,0]$};

                \node[ellipse, line width=1mm, inner sep=0pt,
                    minimum height = 12mm, minimum width = 24mm,
                    draw=black!30] (ll_1) at (-6.5, -2) {$[2,1,2,2,1,0,0]$};
                
                \node[ellipse, line width=1mm, inner sep=0pt,
                    minimum height = 12mm, minimum width = 24mm,
                    draw=black!30] (l_1) at (-2.5, -2) {$[2,2,2,2,1,0,0]$};
                
                \node[ellipse, line width=1mm, inner sep=0pt,
                    minimum height = 12mm, minimum width = 24mm,
                    draw=black!30] (r_1) at (2.5, -2) {$[2,0,2,2,1,0,1]$};
                    
                \node[ellipse, line width=1mm, inner sep=0pt,
                    minimum height = 12mm, minimum width = 24mm,
                    draw=black!30] (rr_1) at (6.5, -2) {$[2,0,2,2,1,0,2]$};

                \node[ellipse, line width=1mm, inner sep=0pt,
                    minimum height = 12mm, minimum width = 24mm,
                    draw=black!30] (lll_2) at (-12, -4) {$[2,1,2,2,1,0,1]$};
                
                \node[ellipse, line width=1mm, inner sep=0pt,
                    minimum height = 12mm, minimum width = 24mm,
                    draw=black!30] (ll_2) at (-8, -4) {$[2,1,2,2,1,0,2]$};
                    
                \node[ellipse, line width=1mm, inner sep=0pt,
                    minimum height = 12mm, minimum width = 24mm,
                    draw=black!30] (rr_2) at (5, -4) {$[2,2,2,2,1,0,1]$};
                
                \node[ellipse, line width=1mm, inner sep=0pt,
                    minimum height = 12mm, minimum width = 24mm,
                    draw=black!30] (rrr_2) at (9, -4) {$[2,2,2,2,1,0,2]$};

                \draw[->, shorten >=1pt] (root) -- (ll_1) node[draw=none,fill=none,midway,below, shift={(1.5,0.3)}] {$a_2$};
                \draw[->, shorten >=1pt] (root) -- (l_1);
                \draw[->, shorten >=1pt] (root) -- (r_1);
                \draw[->, shorten >=1pt] (root) -- (rr_1) node[draw=none,fill=none,midway,below, shift={(-1.5,0.3)}] {$a_7$};
                
                \draw[->, shorten >=1pt] (ll_1) -- (lll_2) node[draw=none,fill=none,midway,right, shift={(1.0, 0.0)}] {$a_7$};
                \draw[->, shorten >=1pt] (ll_1) -- (ll_2);
                \draw[->, shorten >=1pt] (l_1) -- (rr_2);
                \draw[->, shorten >=1pt] (l_1) -- (rrr_2) node[draw=none,fill=none,midway,below] {$a_7$};
                
        \end{tikzpicture}
    }
\end{center}

Recall that the optimal action set $\mathcal{R}^*([2,0,2,2,1,0,0]) = \{a_2, a_7\}$. We must select one of these actions to report in the optimal \gls{dt} during the tree selection step. Let us define our secondary objective function to find the \gls{dt} with the minimal number of nodes. Therefore, the following table summarizes the number of nodes in the optimal subtree of a given state (including itself):

\begin{center}
\begin{tabular}{|c|c c c|}
\hline
     State & Optimal $\mathcal{R}^*(\bs)$ & Nodes in Subtree & Optimal Action\\
     \hline
     $[2,0,2,2,1,0,0]$ & $\{a_2, a_7\}$ & $\min(7, 3) = 3$ & $a_7$\\
     \hline
     $[2,1,2,2,1,0,0]$ & $\{a_7\}$ & $\min(3) = 3$ & $a_7$\\
     $[2,2,2,2,1,0,0]$ & $\{a_7\}$ & $\min(3) = 3$ & $a_7$\\
     $[2,0,2,2,1,0,1]$ & $\emptyset$ & 1 & n/a\\
     $[2,0,2,2,1,0,2]$ & $\emptyset$ & 1 & n/a\\
     \hline
     $[2,1,2,2,1,0,1]$ & $\emptyset$ & 1 & n/a\\
     $[2,1,2,2,1,0,2]$ & $\emptyset$ & 1 & n/a\\
     $[2,2,2,2,1,0,1]$ & $\emptyset$ & 1 & n/a\\
     $[2,2,2,2,1,0,2]$ & $\emptyset$ & 1 & n/a\\
     \hline
\end{tabular}
\end{center}
Given the optimal actions above, the optimal subtree from the state is simply the right branch involving optimal action $a_7$. 

\begin{center}
    \resizebox*{0.5\textwidth}{!}{
        \begin{tikzpicture}
                \node[ellipse, line width=1mm, inner sep=0pt,
                    minimum height = 12mm, minimum width = 24mm,
                    draw=black!30] (root) at (0, 0) {$[2,0,2,2,1,0,0]$};
                
                \node[ellipse, line width=1mm, inner sep=0pt,
                    minimum height = 12mm, minimum width = 24mm,
                    draw=black!30] (r_1) at (-2, -2) {$[2,0,2,2,1,0,1]$};
                    
                \node[ellipse, line width=1mm, inner sep=0pt,
                    minimum height = 12mm, minimum width = 24mm,
                    draw=black!30] (rr_1) at (2, -2) {$[2,0,2,2,1,0,2]$};
                    
                \draw[->, shorten >=1pt] (root) -- (r_1);
                \draw[->, shorten >=1pt] (root) -- (rr_1) node[draw=none,fill=none,midway,below, shift={(-1,0.3)}] {$a_7$};
        \end{tikzpicture}
    }
\end{center}

\end{tcolorbox}

\section{Computational Results}
\label{sec:computationalresults}

\begin{figure}[h!]
    \centering
    \includegraphics[trim={1.75cm 0.5cm 0.75cm 0.5cm},clip, scale=0.75]{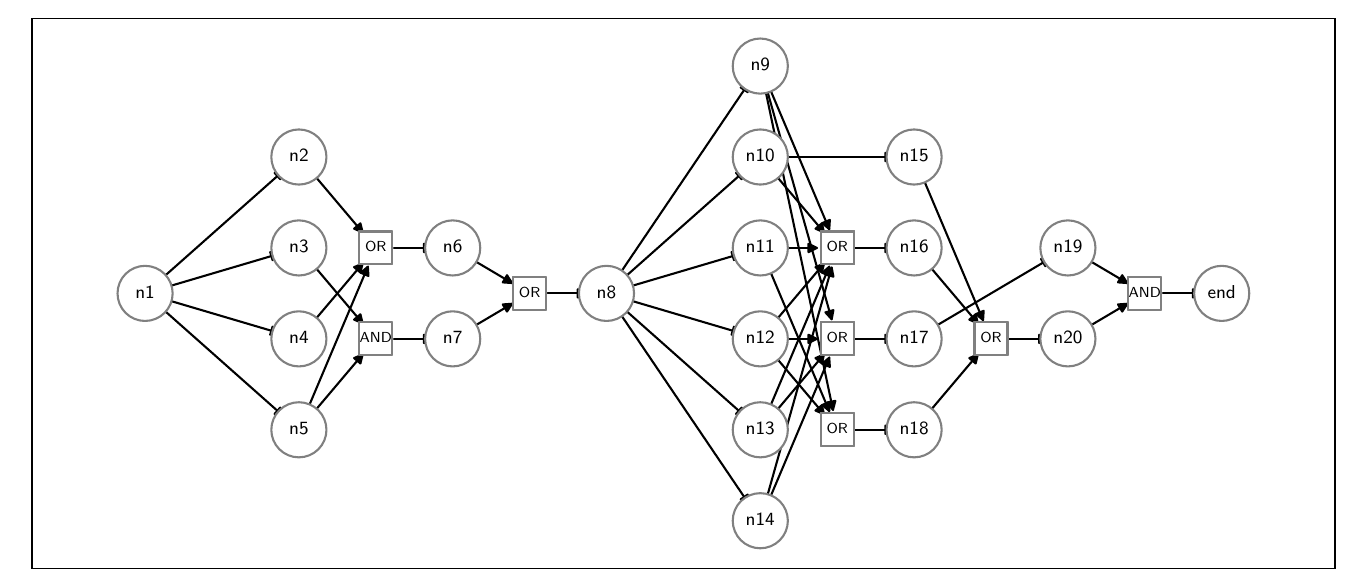}
    \caption{Example random action dependency graph with 20 actions, and OR and AND prerequisite relationships.}
    \label{fig:example_graph}
\end{figure}

In this section, we benchmark our method using randomly generated actions using graph-based AND and OR relationships, similar to the illustrative example as defined in Section~\ref{sec:example}. An example of such a dependency graph shown in Figure~\ref{fig:example_graph}. As a reminder, the circle nodes in the graph represent different actions that the agent may take to get closer to the end node, where rewards are obtained. Different actions have different numbers of outcomes, each of which may enable one or more downstream actions through AND or OR prerequisites, which are highlighted with the square nodes.

In the graph generation process, we make the following simplifying assumptions. We consider prerequisites and but not preclusions, and that each node has either an AND or an OR dependency, but not both. Additionally, we assume that the dependency graph terminates in a single reward node. These assumptions help us generate random scenarios that do not result in infeasible or trivial missions with high probability. We also restrict the length of possible action sequences with an action budget, which sets an absolute limit to the depth of possible \gls{dt}s.

All of the benchmarking has been performed on a personal computer with 32GB of RAM, on a single thread of an Intel i7-11800H processor, using an object-oriented software implementation in Python. Memory was not a limiter in solving any of the benchmarks considered. We make no claims based on the relative efficiency of our implementation; more capable software engineers will likely be able to achieve higher performance through the use of another programming language and more sophisticated data structures. It should also be noted that any variation in benchmarking time within 10\% should be considered to be noise, since there can be substantial variation in the measured time of different benchmarks. In Table \ref{tab:csvbenchmarks}, we provide the computational results over 11 randomized examples, where all times are given in seconds. The columns in the table are as follows: 
\begin{itemize}
    \item Index: Problem specific index, for reference,
    \item $N$: Number of actions,
    \item $B$: Total budget,
    \item $\Phi(\bs_{\rm{root}})$: Reward expectation of the optimal \gls{dt},
    \item $T_{\mathbb{P}}$: Time taken to compute rewarding sets,
    \item $T_{\rm{FG}}$, $T_{\rm{RG}}$, $T_{\rm{DT}}$: Time taken to construct the \gls{fg}, \gls{rg} and optimal \gls{dt}, respectively,
    \item $T_{\rm{total}}$: Total time, 
    \item  $|\mathcal{S}_{\rm{FG}}|$, $|\mathcal{S}_{\rm{RG}}|$,  $|\mathcal{S}_{\rm{DT}}|$: Number of states in the resulting \gls{fg}, \gls{rg} and optimal \gls{dt}, respectively,
    \item $|\mathbb{P}(\bs_{\rm{root}})|$: Total number of rewarding sets. 
\end{itemize}

\begin{table}[ht!]
    \begin{center}
        \resizebox*{\textwidth}{!}{
\begin{tabular}{lrrrrrrrrrrrr}
\toprule
Index & $N$ & $B$ & $\Phi(\bs_{\rm{root}})$ & $T_{\mathbb{P}}$ & $T_{\rm{FG}}$ & $T_{\rm{RG}}$ & $T_{\rm{DT}}$ & $T_{\rm{total}}$ & $|\mathcal{S}_{\rm{FG}}|$ & $|\mathcal{S}_{\rm{RG}}|$ & $|\mathcal{S}_{\rm{DT}}|$ & $|\mathbb{P}(\bs_{\rm{root}})|$ \\
\midrule
0 & 25 & 15 & 0.000243 & 2.09 & 3236.23 & 359.97 & 432.94 & 4031.24 & 2518548 & 2239364 & 235 & 20 \\
1 & 20 & 15 & 0.000971 & 3.00 & 2270.99 & 272.66 & 308.83 & 2855.48 & 1954163 & 1820303 & 309 & 28 \\
2 & 20 & 15 & 0.038782 & 1.70 & 210.13 & 17.65 & 17.45 & 246.94 & 164758 & 105235 & 1347 & 13 \\
3 & 20 & 15 & 0.034325 & 1.62 & 109.27 & 10.00 & 12.82 & 133.72 & 87950 & 76747 & 954 & 12 \\
4 & 15 & 15 & 0.467120 & 1.57 & 14.28 & 1.87 & 2.52 & 20.24 & 13892 & 13892 & 324 & 12 \\
5 & 15 & 10 & 0.466594 & 1.50 & 7.49 & 0.85 & 1.11 & 10.95 & 8186 & 6758 & 314 & 12 \\
6 & 20 & 10 & 0.000968 & 2.87 & 5.06 & 0.34 & 0.17 & 8.43 & 8416 & 1312 & 86 & 27 \\
7 & 17 & 10 & 0.181530 & 1.42 & 5.43 & 0.42 & 0.44 & 7.72 & 8599 & 3188 & 456 & 12 \\
8 & 20 & 10 & 0.026389 & 1.73 & 1.41 & 0.08 & 0.05 & 3.27 & 2237 & 217 & 76 & 12 \\
9 & 25 & 10 & 0.000227 & 2.26 & 0.49 & 0.03 & 0.03 & 2.81 & 1011 & 102 & 40 & 20 \\
10 & 20 & 10 & 0.021860 & 1.37 & 0.85 & 0.06 & 0.08 & 2.36 & 1673 & 344 & 90 & 10 \\
\bottomrule
\end{tabular}
}
\end{center}
    \caption{Benchmarks of networks, with timing in seconds}
    \label{tab:csvbenchmarks}
\end{table}

Experiments in Table \ref{tab:csvbenchmarks} show that the bulk of computation time is dedicated to creating the full graph for the larger examples, i.e. exploring states on the way to rewards. In general, the method breaks the exponential relationship between number of actions and the number of states that need to be explored. This is seen by the weak relationship between $N$ and $|\mathcal{S}_{\rm{FG}}|$. 

Rewarding set computation time $T_{\mathbb{P}}$ is a small proportion of the total time for the more challenging decision problems, and scales linearly with the number of rewarding sets $|\mathbb{P}(\bs_{\rm{root}})|$. However, it should be noted that, for the smaller benchmarks, this computation tends to be the slowest component of the algorithm, since it relies on solving an inefficient \gls{mio} with cutting planes as defined in Appendix~\ref{sec:rewardingsetsgraph}. A constraint satisfaction problem approach would likely outperform the \gls{mio} we use in this paper, though the marginal gains are small in the context of total time.

To give an idea of the magnitude of speed-ups due to the use of rewarding sets, we have included the time taken to generate the optimal solution using the naive approach in Algorithm~\ref{alg:Naive_FG} for a subset of the examples in Table~\ref{tab:csvbenchmarksnaive}. This subset includes problems that could be solved within 4000 seconds, which is roughly the maximum time it took to solve any benchmark using the accelerated methodology. Note that for the accelerated algorithm, the time overhead for the computation of rewarding sets is included in the total time, in order to have a fair comparison. 

\begin{table}[ht!]
    \centering
        \resizebox*{0.8\textwidth}{!}{
\begin{tabular}{lrrrrrrrrrrr}
\toprule
Index & $N$ & $B$ & $\Phi(\bs_{\rm{root}})$ & $T_{\rm{FG}}$ & $T_{\rm{RG}}$ & $T_{\rm{DT}}$ & $T_{\rm{total}}$ & $|\mathcal{S}_{\rm{FG}}|$ & $|\mathcal{S}_{\rm{RG}}|$ & $|\mathcal{S}_{\rm{DT}}|$ \\
\midrule
5 & 15 & 10 & 0.466594 & 7.24 & 0.76 & 0.95 & 8.95 & 9466 & 6758 & 314 \\
7 & 17 & 10 & 0.181530 & 21.53 & 1.49 & 0.55 & 23.57 & 18775 & 4064 & 456 \\
8 & 20 & 10 & 0.026389 & 155.15 & 7.98 & 0.08 & 163.21 & 83854 & 525 & 76 \\
10 & 20 & 10 & 0.021860 & 78.98 & 4.53 & 0.11 & 83.62 & 47320 & 760 & 90  \\
\bottomrule
\end{tabular}
}
    \caption{Benchmarks of networks using the naive approach in  Algorithm~\ref{alg:Naive_FG}.}
    \label{tab:csvbenchmarksnaive}
\end{table}

Solutions from~Table~\ref{tab:csvbenchmarksnaive} match those from Table~\ref{tab:csvbenchmarks} in both $\Phi(\bs_{\rm{root}})$ and $|\mathcal{S}_{\rm{FG}}|$, giving us confidence that both methods find the same optimal \gls{dt}s. Though the size of the \gls{rg} is not consistent between methods, this is expected behavior. Non-rewarding actions can be present in the \gls{rg} when reward-based pruning is not used, since these actions do not impact optimality. These extraneous actions are removed when determining the optimal \gls{dt}, since the default secondary objective is to minimize tree size as described in Section~\ref{sec:optimaltree}.

With the naive approach, we see the expected increasing relationship between the number of nodes and the size of the full graph, and thus computational time. The maximum speed-up observed from the naive approach to the accelerated approach is 50 times, in the solution of problem 10.

\begin{figure}[h!]
    \centering
    \includegraphics[trim={0.0cm 0cm 1.0cm 1.0cm}, clip, scale=0.9]{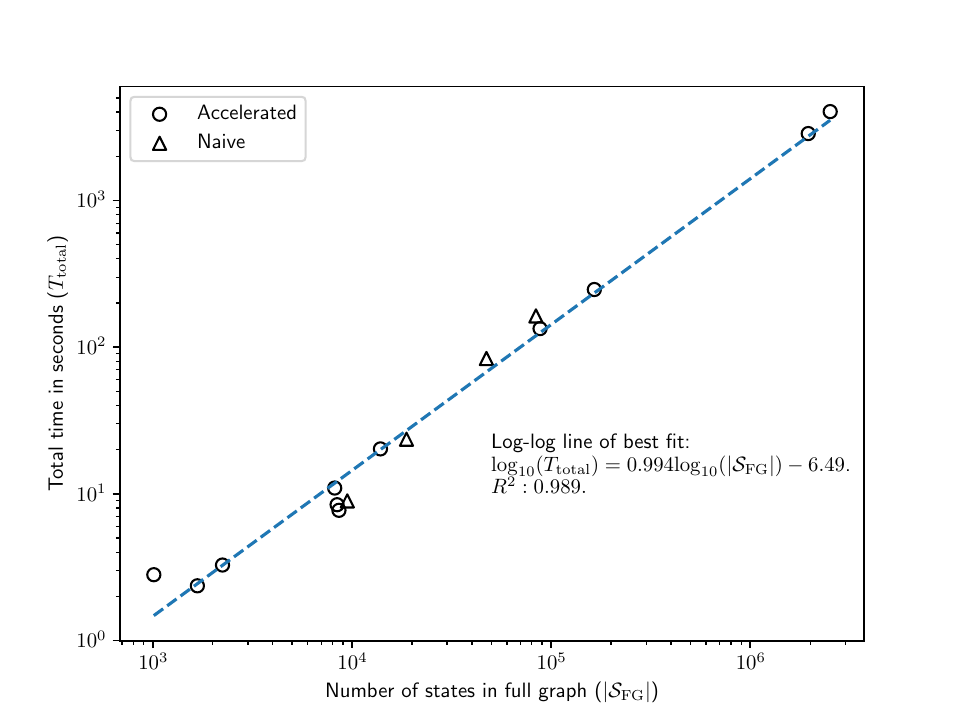}
    \caption{The optimal \gls{dt} method scales linearly with the number of states explored.}
    \label{fig:scaling}
\end{figure}

An interesting case study is problem 5, where the naive approach outperforms our method, with total times of 8.95 to 10.95 seconds, respectively. This is due to the fact that, for this problem, the action pruning method doesn't reduce the size of the state space substantially, but adds overhead due to the time required to compute the rewarding sets and propagate them through the \gls{fg}. Problem 5 is a special case that demonstrates that, while the approach provides substantial speed-ups for a majority of decision problems, there can be cases where no or negative benefit is observed. In general, bigger improvements are expected as the size of the decision problem grows, as the number of pruned states grows larger and the upfront computational time for computing rewarding sets takes a diminishing proportion of the total time.

When trying to estimate the total time, we observe that the driving factor is not the number of actions. The dependencies between actions make it difficult to predict a priori the number of strategies available to two different problems with the same number of actions but different \gls{adg} structure. Thus, no direct relationship is observed between $N$, $|\mathbb{P}(\bs_{\rm{root}})|$ and $T_{\rm{total}}$. Instead, the number of states in the full graph drives the computation time. Figure~\ref{fig:scaling} shows that total time scales linearly with the number of states explored, where the log-log line of best fit has been computed for the accelerated algorithm, with an $R^2$ of 0.989. The naive approach shows similar scaling behavior, and its benchmarks have been included for comparison.

\section{Conclusion}

We have introduced a methodology to generate optimal \gls{dt}s to make decisions in settings with interdependent actions with probabilistic outcomes. The methodology leverages dynamic programming and mixed-integer linear optimization methods, and uses problem-specific information to find \gls{dt}s efficiently without sacrificing global optimality. Through computational experiments, we showed that this method's computation time scales linearly with the number of states explored. Therefore, a natural extension of this research is to further reduce the number of states required to explore while maintaining optimality. 

Additionally, we have also developed some nascent techniques for the \emph{conditioning} of actions (i.e. previous actions and outcomes changing the parameters of another action) which is amenable to the algorithms in Section~\ref{sec:methodology} and could extend the applications of \gls{dt}s. These techniques would extend the existing method and this research is heavily focused on the efficient implementation of the idea. Finally, we are actively working on an adversarial variant of these algorithms, both designing \gls{coa}s in the presence of an adversary or uncertainty and also determining the optimal strategy to protect against an optimal \gls{dt}.

\section{Acknowledgments}
\label{sec:acknowledgements}

The authors thank Mr. John Garstka, Director for Cyber Warfare at the Office of Department of Defense (Acquisition and Sustainment) for his vision in developing the Cyber Warfare Analysis domain and his support for this paper, as well as Dr.\ Doug Altner (MITRE) for his feedback. 

\bibliographystyle{apalike}
\bibliography{main.bib}
\section{Appendix}

\subsection{Rewarding set computation for action dependency graphs}
\label{sec:rewardingsetsgraph}

In order to generalize the rewarding set computation for an \gls{adg} with preclusions and prerequisites, we have to introduce the following new machinery:
\begin{itemize}
    \item Function $\rm{AND}: r \xrightarrow{} \omega \subseteq \Omega$ and function $\rm{OR}: r \xrightarrow{} \omega \subseteq \Omega$ describe the sets of actions and outcomes that are and-prerequisites and or-prerequisites of the action $r$,
    \item Function $\rm{NOTAND}: r \xrightarrow{} \omega \subseteq \Omega$ and function $\rm{NOTOR}: r \xrightarrow{} \omega \subseteq \Omega$ describe the sets of actions and outcomes that are not-and-prerequisites and not-or-prerequisites of the action $r$, and
    \item $\Omega_{\rho} \subset \Omega$ describes the action-outcome pairs that immediately precede a reward in the \gls{adg}\footnote{If more complicated relationships define the reward, it is trivial to add a rewarding dummy action with a single outcome and a cost of $0$ with the necessary prerequisites to define the relationship.}. In the illustrative example from Section~\ref{sec:example}, $\Omega_{\rho} = \{(a_5, 2), (a_6, 2), (a_7, 2)\}$. 
\end{itemize}

While the non-dominated reward states $\mathcal{S}_{\rho}$ are hard to identify directly, each rewarding set $\mathcal{P}(\bs_{\rm{root}}, \bs_k)$ can be identified, one dominating rewarding state $\bs_k$ at a time, by iteratively solving a shortest path algorithm on the \gls{adg} to find action-outcome trajectories from the start to the reward nodes. Since this has to be done while considering the various action prerequisites and preclusions, it can only be expressed via the following \gls{mio} problem over all possible action and outcome combinations $\Omega$:
\begin{align}
    \underset{\bz}{\rm{minimize}} \quad & \be^T\bz \label{eq:MIOPathobj}\\ 
    \rm{subject~to~}\quad & \sum_{(r,j) \in \Omega_{\rho}} z_{(r,j)} \geq 1 \label{eq:sinks} \\ 
    & \sum_{(r_i, j_i)~ \in ~\rm{OR}(r_o, j_o)} z_{(r_i,j_i)} \geq z_{(r_o, j_o)},~\forall  (r_o, j_o) \in \Omega,\\
    & \sum_{(r_i, j_i)~ \in ~\rm{NOTOR}(r_o, j_o)} (1-z_{(r_i,j_i)}) \geq z_{(r_o, j_o)},~\forall  (r_o, j_o) \in \Omega,\\
    & z_{(r_i,j_i)} \geq z_{(r_o, j_o)}, ~\forall (r_o, j_o) \in \Omega, (r_i, j_i) \in \rm{AND}(r_o, j_o),\\
    & (1 - z_{(r_i,j_i)}) \geq z_{(r_o, j_o)}, ~\forall (r_o, j_o) \in \Omega, (r_i, j_i) \in \rm{NOTAND}(r_o, j_o),\\
    & z_{(r,j)} \in \{0, 1\}, (r,j) \in \Omega. \label{eq:z}
\end{align}
The optimal solution to the above finds one rewarding set through the value of $\bz^*$, given by 
\begin{equation}
    \mathcal{P}(\bs_{\rm{root}}, \bs_k) = \{(r,j) \subseteq \Omega: z^*_{(r,j)} = 1\}
    \label{eq:bzk}
\end{equation}
from which the individual $\bs_k \in \mathcal{S}_{\rho}$ can be extracted via the rewarding set to rewarding state relationship in Equation~\eqref{eq:pathrelation}. 

To find the next $\mathcal{P}(\bs_{\rm{root}}, \bs_k)$ with $\bs_k \in \mathcal{S}_{\rho}$, we must ensure that the rewarding states from previous iterations do not dominate $\bs_k$. To do so, we add a cut to exclude the previous solution $\bz^*$ from the optimization problem,
\begin{equation}
    \sum_{\Big\{(r,j) \subseteq \Omega :~z^*_{(r,j)} = 1\Big\}} z_{(r,j)} \leq |\bz^*| - 1, 
\end{equation}
and re-optimize until the \gls{mio} becomes infeasible. The set containing all $\mathcal{P}(\bs_{\rm{root}}, \bs_k)$ computed using Equation~\eqref{eq:bzk} over the \gls{mio} instances is $\mathbb{P}(\bs_{\rm{root}})$, and the set containing the corresponding $\bs_k$ is the the set of dominating reward states $\mathcal{S}_{\rho}$.

\subsection{State indexing}
\label{sec:stateindexing}
One challenge when dealing with a large number of states is to be able to store and query them efficiently. We have developed a method to generate a unique numerical index for each state based on the domain of each element of the state vector, 
\begin{equation}
    \rm{dom}(\mathbb{S}_i),~\forall i \in n,
\end{equation}
extending the concept of domain for each dimension of $\mathbb{S} \in \mathbb{Z}^n$. Given that each $\rm{dom}(\mathbb{S}_i)$ is finite, we can uniquely identify each possible state through the following formula for state index, 
\begin{equation} 
    \rm{index}(\bs) = \sum_{i=1}^{n} \Bigg(\sum_{x \in \rm{dom}(\mathbb{S}_i)} \mathbb{I}(x \leq s_i) \Bigg) \prod_{j=1}^{i-1} |\rm{dom}(\mathbb{S}_j)|.
\end{equation} 
This is equivalent to creating a custom basis. As a concrete example, if the possible states were belonging to a vector in $\{0, 1, \ldots, 9\}^5$, this would be equivalent to mapping each state to an integer value between 0 and 99999. However, since in practice different elements of the state have different bases, we create a custom basis that is unique to any valid instance of the state vector.

This allows us to store states in our graphs or trees efficiently in a dictionary. In general, the big-$\mathcal{O}$ time taken to query a well-hashed dictionary is $\mathcal{O}(1)$ with the worst case being  $\mathcal{O}(n)$, where $n$ is the number of elements in the dictionary. For comparison, for a list, the average case has the same performance as the worst case, which is $\mathcal{O}(n)$. 

\end{document}